\documentclass[russian,a4paper,11pt]{article}

\usepackage[russian]{babel}

\usepackage[all]{xy}
\usepackage[dvips]{graphicx}
\usepackage{amsmath,amssymb,amsfonts,amsthm,latexsym,mathrsfs,textcomp,verbatim, enumerate, soul, color}
\xyoption{web}

\title{Унификация математики\\ с помощью теории топосов}
\author{Оливия Карамелло\thanks{Supported by a visiting position from the Centro di Ricerca Matematica Ennio De Giorgi, Pisa.}}

\date{20 июня 2010 г.}

\begin{document}

\mathcode`\<="4268  
\mathcode`\>="5269  
\mathcode`\.="313A  
\mathchardef\semicolon="603B 
\mathchardef\gt="313E
\mathchardef\lt="313C

\newcommand{\app}
 {{\sf app}}

\newcommand{\Ass}
 {{\bf Ass}}

\newcommand{\ASS}
 {{\mathbb A}{\sf ss}}

\newcommand{\Bb}
{\mathbb}

\newcommand{\biimp}
 {\!\Leftrightarrow\!}

\newcommand{\bim}
 {\rightarrowtail\kern-1em\twoheadrightarrow}

\newcommand{\bjg}
 {\mathrel{{\dashv}\,{\vdash}}}

\newcommand{\bstp}[3]
 {\mbox{$#1\! : #2 \bim #3$}}

\newcommand{\cat}
 {\!\mbox{\t{\ }}}

\newcommand{\cinf}
 {C^{\infty}}

\newcommand{\cinfrg}
 {\cinf\hy{\bf Rng}}

\newcommand{\cocomma}[2]
 {\mbox{$(#1\!\uparrow\!#2)$}}

\newcommand{\cod}
 {{\rm cod}}

\newcommand{\comma}[2]
 {\mbox{$(#1\!\downarrow\!#2)$}}

\newcommand{\comp}
 {\circ}

\newcommand{\cons}
 {{\sf cons}}

\newcommand{\Cont}
 {{\bf Cont}}

\newcommand{\ContE}
 {{\bf Cont}_{\cal E}}

\newcommand{\ContS}
 {{\bf Cont}_{\cal S}}

\newcommand{\cover}
 {-\!\!\triangleright\,}

\newcommand{\cstp}[3]
 {\mbox{$#1\! : #2 \cover #3$}}

\newcommand{\Dec}
 {{\rm Dec}}

\newcommand{\DEC}
 {{\mathbb D}{\sf ec}}

\newcommand{\den}[1]
 {[\![#1]\!]}

\newcommand{\Desc}
 {{\bf Desc}}

\newcommand{\dom}
 {{\rm dom}}

\newcommand{\Eff}
 {{\bf Eff}}

\newcommand{\EFF}
 {{\mathbb E}{\sf ff}}

\newcommand{\empstg}
 {[\,]}

\newcommand{\epi}
 {\twoheadrightarrow}

\newcommand{\estp}[3]
 {\mbox{$#1 \! : #2 \epi #3$}}

\newcommand{\ev}
 {{\rm ev}}

\newcommand{\Ext}
 {{\rm Ext}}

\newcommand{\fr}
 {\sf}

\newcommand{\fst}
 {{\sf fst}}

\newcommand{\fun}[2]
 {\mbox{$[#1\!\to\!#2]$}}

\newcommand{\funs}[2]
 {[#1\!\to\!#2]}

\newcommand{\Gl}
 {{\bf Gl}}

\newcommand{\hash}
 {\,\#\,}

\newcommand{\hy}
 {\mbox{-}}

\newcommand{\im}
 {{\rm im}}

\newcommand{\imp}
 {\!\Rightarrow\!}

\newcommand{\Ind}[1]
 {{\rm Ind}\hy #1}

\newcommand{\iten}[1]
{\item[{\rm (#1)}]}

\newcommand{\iter}
 {{\sf iter}}

\newcommand{\Kalg}
 {K\hy{\bf Alg}}

\newcommand{\llim}
 {{\mbox{$\lower.95ex\hbox{{\rm lim}}$}\atop{\scriptstyle
{\leftarrow}}}{}}

\newcommand{\llimd}
 {\lower0.37ex\hbox{$\pile{\lim \\ {\scriptstyle
\leftarrow}}$}{}}

\newcommand{\Mf}
 {{\bf Mf}}

\newcommand{\Mod}
 {{\bf Mod}}

\newcommand{\MOD}
{{\mathbb M}{\sf od}}

\newcommand{\mono}
 {\rightarrowtail}

\newcommand{\mor}
 {{\rm mor}}

\newcommand{\mstp}[3]
 {\mbox{$#1\! : #2 \mono #3$}}

\newcommand{\Mu}
 {{\rm M}}

\newcommand{\name}[1]
 {\mbox{$\ulcorner #1 \urcorner$}}

\newcommand{\names}[1]
 {\mbox{$\ulcorner$} #1 \mbox{$\urcorner$}}

\newcommand{\nml}
 {\triangleleft}

\newcommand{\ob}
 {{\rm ob}}

\newcommand{\op}
 {^{\rm op}}

\newcommand{\palrr}[4]{
  \def\labelstyle{\scriptstyle}
  \xymatrix{ {#1} \ar@<0.5ex>[r]^{#2} \ar@<-0.5ex>[r]_{#3} & {#4} } }

\newcommand{\palrl}[4]{
  \def\labelstyle{\scriptstyle}
  \xymatrix{ {#1} \ar@<0.5ex>[r]^{#2}  &  \ar@<0.5ex>[l]^{#3} {#4} } }

\newcommand{\pepi}
 {\rightharpoondown\kern-0.9em\rightharpoondown}

\newcommand{\pmap}
 {\rightharpoondown}

\newcommand{\Pos}
 {{\bf Pos}}

\newcommand{\prarr}
 {\rightrightarrows}

\newcommand{\princfil}[1]
 {\mbox{$\uparrow\!(#1)$}}

\newcommand{\princid}[1]
 {\mbox{$\downarrow\!(#1)$}}

\newcommand{\prstp}[3]
 {\mbox{$#1\! : #2 \prarr #3$}}

\newcommand{\pstp}[3]
 {\mbox{$#1\! : #2 \pmap #3$}}

\newcommand{\relarr}
 {\looparrowright}

\newcommand{\rlim}
 {{\mbox{$\lower.95ex\hbox{{\rm lim}}$}\atop{\scriptstyle
{\rightarrow}}}{}}

\newcommand{\rlimd}
 {\lower0.37ex\hbox{$\pile{\lim \\ {\scriptstyle
\rightarrow}}$}{}}

\newcommand{\rstp}[3]
 {\mbox{$#1\! : #2 \relarr #3$}}

\newcommand{\scn}
 {{\bf scn}}

\newcommand{\scnS}
 {{\bf scn}_{\cal S}}

\newcommand{\semid}
 {\rtimes}

\newcommand{\Sep}
 {{\bf Sep}}

\newcommand{\sep}
 {{\bf sep}}

\newcommand{\Set}
 {{\bf Set}}

\newcommand{\Sh}
 {{\bf Sh}}

\newcommand{\ShE}
 {{\bf Sh}_{\cal E}}

\newcommand{\ShS}
 {{\bf Sh}_{\cal S}}

\newcommand{\Simp}
 {{\bf \Delta}}

\newcommand{\snd}
 {{\sf snd}}

\newcommand{\stg}[1]
 {\vec{#1}}

\newcommand{\stp}[3]
 {\mbox{$#1\! : #2 \to #3$}}

\newcommand{\Sub}
 {{\rm Sub}}

\newcommand{\SUB}
 {{\mathbb S}{\sf ub}}

\newcommand{\tbel}
 {\prec\!\prec}

\newcommand{\tic}[2]
 {\mbox{$#1\!.\!#2$}}

\newcommand{\tp}
 {\!:}

\newcommand{\tps}
 {:}

\newcommand{\tsub}
 {\pile{\lower0.5ex\hbox{.} \\ -}}

\newcommand{\wavy}
 {\leadsto}

\newcommand{\wavydown}
 {\,{\mbox{\raise.2ex\hbox{\hbox{$\wr$}
\kern-.73em{\lower.5ex\hbox{$\scriptstyle{\vee}$}}}}}\,}

\newcommand{\wbel}
 {\lt\!\lt}

\newcommand{\wstp}[3]
 {\mbox{$#1\!: #2 \wavy #3$}}

\newcommand{\fu}[2]
{[#1,#2]}


%
%
%
\def\pushright#1{{
   \parfillskip=0pt            
   \widowpenalty=10000         
   \displaywidowpenalty=10000  
   \finalhyphendemerits=0      
  %
   \leavevmode                 
   \unskip                     
   \nobreak                    
   \hfil                       
   \penalty50                  
   \hskip.2em                  
   \null                       
   \hfill                      
   {#1}                        
  %
   \par}}                      

\def\qed{\pushright{$\square$}\penalty-700 \smallskip}

\newtheorem{theorem}{Theorem}[section]

\newtheorem{proposition}[theorem]{Proposition}

\newtheorem{scholium}[theorem]{Scholium}

\newtheorem{lemma}[theorem]{Lemma}

\newtheorem{corollary}[theorem]{Corollary}

\newtheorem{conjecture}[theorem]{Conjecture}

\newenvironment{proofs}%
 {\begin{trivlist}\item[]{\bf Proof }}%
 {\qed\end{trivlist}}

  \newtheorem{rmk}[theorem]{Remark}
\newenvironment{remark}{\begin{rmk}\em}{\end{rmk}}

  \newtheorem{rmks}[theorem]{Remarks}
\newenvironment{remarks}{\begin{rmks}\em}{\end{rmks}}

  \newtheorem{defn}[theorem]{Definition}
\newenvironment{definition}{\begin{defn}\em}{\end{defn}}

  \newtheorem{eg}[theorem]{Example}
\newenvironment{example}{\begin{eg}\em}{\end{eg}}

  \newtheorem{egs}[theorem]{Examples}
\newenvironment{examples}{\begin{egs}\em}{\end{egs}}


\bgroup           
\let\footnoterule\relax  
\maketitle

\begin{center}
 {\scriptsize Перевод с английского\\}
 {\footnotesize Ю.~Н.~Браткова}
\end{center}

\vspace{0.25cm}

\begin{abstract}
Дан набор принципов и методологий, могущих служить основанием
унифицирующей теории математики. Эти принципы основаны на новом
видении топосов Гротендика как унифицирующих пространств,
способных служить 'мостами' для переноса информации, идей и
результатов между различными математическими теориями.
\end{abstract}
\egroup

\vspace{4.5cm}
{\small \textsc{Centro di Ricerca Matematica Ennio De Giorgi, Scuola Normale Superiore, Piazza dei Cavalieri 3, 56100 Pisa, Italy}\\
\emph{E-mail address:} \texttt{olivia.caramello@sns.it}}

\vspace{0.5cm}

{\small \textsc{Department of Pure Mathematics and Mathematical Statistics, University of Cambridge, Wilberforce Road, Cambridge CB3 0WB, UK}\\
\emph{E-mail address:} \texttt{O.Caramello@dpmms.cam.ac.uk}}

\vspace{5 mm}

\tableofcontents

\newpage

\section{Введение}

Цель этой статьи --- установить набор принципов и методологий,
могущих служить основанием унифицирующей теории математики. Под
этим мы подразумеваем метатеорию, предлагающую методологии для
сравнения различных математических теорий друг с другом, как для
раскрытия аналогии между ними, так и для выделения особенностей,
и, что наиболее важно, дающую эффективные смыслы для переноса
результатов и техники между различными областями.

Математика подразделяется на несколько различных областей:
алгебра, анализ, геометрия, топология, теория чисел и т.д. Каждая
из этих областей развивается многие годы, разрабатывая свои
собственные идеи и технику, и достигла к нынешнему моменту высокой
степени специализации. Со временем открываются различные связи
между областями, что приводит в некоторых случаях к созданию
настоящих `мостов' между разными ветвями математики; множество раз
методы одной области применялись для получения результатов в
другой, и это взаимодействие в одной и той же теме различных точек
зрения всегда играло фундаментальную роль в выявлении природы
концепций, предлагая новые результаты и подсказывая новые линии
исследования.

Глубокие причины, стоящие за аналогиями и связями между различными
областями, лучше всего могут быть поняты в контексте
математической логики; в самом деле, логика дает осмысление для
формализаций любой разновидности математических концепций, поэтому
исследование взаимоотношений между различными теориями может быть
вынесено на полностью строгий уровень.

С точки зрения важности наведения мостов между разными ветвями
математики было бы крайне желательно, чтобы логика служила бы не
столько инструментом анализа уже обнаруженных в математике
аналогий, сколько играла бы активную роль в установлении новых
связей между существующими областями, а также подсказывала бы
новые направления математического исследования. Получилось так,
что сейчас мы имеем математические инструменты для наших намерений
достичь этой цели.

Теория множеств как система, в рамках которой все обычные
математические концепции могут быть строго выражены, представляла
собой первую серьезную попытку логики унифицировать математику по
крайней мере на уровне языка. Позднее теория категорий предложила
альтернативный абстрактный язык, в рамках которого могла быть
сформулирована б\'{о}льшая часть математики, и, таким образом,
стала дальнейшим продвижением в направлении заявленной цели
`унификации математики'. Тем не менее обе эти системы осуществляют
унификацию, весьма ограниченную по своим возможностям в том
смысле, что даже хотя каждая из них дает способ выражения и
организации математики при помощи единого языка, они не предлагают
сами по себе эффективных методов для переноса знаний между
различными областями. С другой стороны, принципы, которые мы
представим в настоящей статье, определяют другой, более
основательный подход к унификации математики.

Наша методология базируется на новом взгляде на топосы Гротендика
как на унифицирующие пространства, которые могут служить `мостами'
для переноса информации, идей и результатов между различными
математическими теориями. Этот взгляд на топосы возник в ходе
диссертационных исследований, проводившихся автором в Кембриджском
университете с 2006 по 2009 гг., и фактически результаты~\cite{OC}
дают вынужденное техническое свидетельство обоснованности такой
точки зрения (в самом деле, одна из целей настоящей статьи ---
служить `концептуальным проводником' к методам авторской
диссертации).

В этой статье я даю набросок фундаментальных принципов,
характеризующих мое видение топосов как `мостов', связывающих
различные математические теории, и описываю общую методологию,
которая возникла из такого видения и которая до сих пор мотивирует
мои исследования. Анализ будет завершен обсуждением
результатов~\cite{OC}, которые лучше всего иллюстрируют применение
упомянутых выше принципов. Эти принципы абстрактны и
трансверсальны к различным математическим областям, и применение
их может привести к огромному количеству удивительных и
нетривиальных результатов в любой области математики, поэтому мы
надеемся, что читатель будет иметь мотивировки для апробации этих
методов в своей области интересов.

\section{Геометрические теории\\ и их классифицирующие топосы}

В этой секции мы обсуждаем основания геометрических теорий и
классифицирующих топосов, необходимые для дальнейших частей статьи
(в качестве исчерпывающего источника информации по теме читателю
рекомендуется \cite{El}, часть~D).

\subsection{Геометрические теории}

\emph{Геометрические теории} --- широкий класс (многосортных)
бесконечных теорий первого порядка. Напомним, что теория над
сигнатурой $\Sigma$ геометрическая, если ее аксиомы могут быть
представлены в форме ${\forall \vec{x}(\phi \imp \psi)}$, где
$\phi$ и $\psi$ --- геометрические формулы над $\Sigma$ (т.е.
формулы с конечным числом свободных переменных, построенные из
атомарных формул над $\Sigma$ с использованием только лишь
конечных конъюнкций, бесконечных дизъюнкций и кванторов
существования в контексте $\vec{x}$. Аксиомы геометрических теорий
часто представлены в последовательной форме, т.е. пишут ${\phi
\vdash_{\vec{x}} \psi}$ вместо ${(\forall \vec{x})(\phi \imp
\psi)}$. Если в вышеупомянутом определении геометрической теории
заменить `бесконечные дизъюнкции' на `конечные дизъюнкции',
получим класс теорий, известных как \emph{когерентные теории}.

Атрибут `геометрический' не означает, что это класс теорий,
которые имеют какое-то особое отношение к геометрии (в отличие от
исторически первых приложений); это вполне общая концепция, и
геометрические теории могут быть обнаружены в любой области
математики.

Действительно, факты таковы, что большинство теорий, естественно
возникающих в математике, имеет геометрическую аксиоматизацию (над
своей сигнатурой). Как бы то ни было, если конечная теория первого
порядка $\mathbb T$ не геометрическая, мы можем канонически
сконструировать когерентную теорию над более широкой сигнатурой,
называемую \emph{морлеизацией} $\mathbb T$, модель которой в
категории множеств $\Set$ (более общо, в любой булевой когерентной
категории) может быть отождествлена с $\mathbb T$ (см. \cite{El},
лемма D1.5.13).

Понятие морлеизации важно, потому что оно позволяет нам изучать
любую теорию первого порядка с использованием методов теории
топосов. По существу мы можем ожидать, что многие свойства теорий
первого порядка естественным образом выражаются как свойства их
морлеизаций, и что эти последние свойства окажутся выразимыми в
терминах `инвариантов' их классифицирующих топосов. Например,
теория первого порядка полна тогда и только тогда, когда ее
морлеизация ${\mathbb T}'$ удовлетворяет следующему свойству:
любое когерентное (эквивалентно, геометрическое) высказывание над
его сигнатурой доказуемо эквивалентно $\bot$ или $\top$, но не
тому и другому; и это свойство точно эквивалентно словам о том,
что классифицирующий топос теории двузначен.

Благодаря своей бесконечной природе геометрическая логика весьма
выразительна; важные математические свойства, которые не
выражаются в конечной логике первого порядка (над данной
сигнатурой), часто допускают геометрическую аксиоматизацию.
Например, свойство элемента коммутативного кольца с единицей быть
нильпотентным, как хорошо известно, не выражается финитарной
формулой первого порядка над сигнатурой колец, но оно с
очевидностью выражается как геометрическая формула над этой
сигнатурой. Интересная бесконечная геометрическая теория, которая
изучалась теоретико-топосно в \cite{OC4}, --- теория полей
конечной характеристики, которые являются алгебраическими над
своими первичными (prime) полями. Другой пример геометрической
теории, которая не является конечной теорией первого порядка, дает
теория абелевых групп с кручением.

Геометрические теории, включая бесконечные, определенно должны
рассматриваться как объекты, значимые для исследования.
Бесконечная природа геометрической логики должна представляться
помехой к распознанию этого класса теорий как субъекта первичной
важности; фактически никто не может, вообще говоря, использовать
методы классической теории моделей для изучения геометрических
теорий, и, следовательно, должны быть использованы методы
совершенно другой природы.

Как мы увидим в разделе~\ref{logGrothendieck}, с точки зрения
теории топосов предположения конечности по отношению к природе
логики действительно представляются весьма неестественными
ограничениями; даже если начать с конечных теорий,
теоретико-топосное их изучение часто приводит к рассмотрению
бесконечных теорий. Логика, лежащая в основе топосов Гротендика,
--- это геометрическая логика в своей полностью бесконечной природе (см.
раздел~\ref{logGrothendieck}); и, как мы покажем в дальнейшем,
теория топосов предоставляет набор невероятно мощной техники для
изучения геометрических теорий.

Подчеркнем, что геометрические теории, как и любая разновидность
теорий первого порядка, --- объекты чисто синтаксической природы.
Как и в классической конечной логике первого порядка, мы имеем
дело с синтаксическим обозначением доказуемости предложений
первого порядка (по отношению к теории), поэтому мы имеем
естественную систему доказательств для геометрической
(соответственно, когерентной) логики, описанной в терминах правил
вывода, использующих геометрические (соответственно, когерентные)
секвенции, которые влекут за собой понятие доказуемости
геометрических (соответственно, когерентных) секвенций по
отношению к данной геометрической (соответственно, когерентной)
теории. В геометрической логике классическая и интуиционистская
доказуемость геометрических секвенций совпадают, поэтому мы вполне
можем опустить закон исключенного третьего из этих систем
доказательства, не затрагивая соответствующее понятие
доказуемости. Для детального представления об этих системах мы
отсылаем читателя к \cite{El}, часть~D.

Хорошо известно, что языки первого порядка всегда могут быть
интерпретированы в контексте (данной модели) теории множеств.
Фактически эти языки могут также быть осмысленно интерпретированы
в категории при условии, что последняя обладает достаточной
категорной структурой, позволяющей интерпретировать связки и
кванторы, возникающие в формулах языка (в этой категорной
семантике типы (sorts) интерпретируются как объекты, термы как
стрелки и формулы как подобъекты в том, что касается логической
структуры составных выражений); это приводит, при тех же самых
предположениях, к понятию удовлетворения секвенций в
категориальной структуре и, следовательно, к понятию модели теории
первого порядка в категории, которая приспособлена, в случае
категории множеств, к классическому определению Тарского
(основанной на множествах) модели теории первого порядка.
Например, топологическая группа может рассматриваться как модель
теории групп в категории топологических пространств, а пучок колец
на топологическом пространстве $X$ может рассматриваться как
модель теории колец в топосе $\Sh(X)$. Эта категориальная
семантика, когда она определена, всегда созвучна вышеупомянутой
системе доказательств при условии, что последняя не содержит
закона исключенного третьего.

Топосы Гротендика --- широкий класс категорий, в которых мы можем
интерпретировать геометрическую логику, и фактически они
обеспечивают нам сильную форму полноты для геометрических теорий.
Действительно, для любой геометрической теории $\mathbb T$ имеется
топос Гротендика $\Set[{\mathbb T}]$, а именно
\emph{классифицирующий топос} для $\mathbb T$, который содержит
\emph{консервативную} модель для $\mathbb T$ (т.е. модель
$U_{\mathbb T}$ такую, что геометрические секвенции, доказуемые в
$\mathbb T$, в точности те же самые, которые удовлетворяются в
$U_{\mathbb T}$), что является, кроме того, \emph{универсальным} в
том смысле, что любая модель для $\mathbb T$ в топосе Гротендика
возникает, с точностью до изоморфизма, как образ объекта
$U_{\mathbb T}$ относительно функтора обратного образа для
(единственного с точностью до эквивалентности) геометрического
морфизма топосов из топоса, в котором модель живет, в
классифицирующий топос (см. раздел \ref{classif} ниже). Отметим,
что, согласно своей бесконечной природе, геометрические теории,
вообще говоря, не обладают классической формой полноты, т.е.,
вообще говоря, неверно, что секвенция, которая действительна во
всех основанных на множествах моделях геометрической теории,
доказуема в теории с использованием геометрической логики. В любом
случае, как было сказано выше, концепция универсальной модели
геометрической теории влечет за собой сильную форму полноты, и
фактически мы верим, что ее хорошо было бы взять в качестве
фундаментального инструмента (замещающего стандартные) в
исследовании аспектов полноты теории (см. разделы \ref{classif} и
\ref{logGrothendieck}).

\subsection{Классифицирующие топосы}\label{classif}

\emph{Классифицирующий топос} геометрической теории $\mathbb T$
над сигнатурой $\Sigma$ --- это топос Гротендика $\Set[{\mathbb
T}]$ такой, что для любого топоса Гротендика $\cal E$ категория
${\bf Geom}({\cal E}, \Set[{\mathbb T}])$ геометрических морфизмов
из ${\cal E}$ в $\Set[{\mathbb T}]$ эквивалентна категории моделей
для $\mathbb T$ в топосе $\cal E$ естественно в $\cal E$;
естественность означает, что для любого геометрического морфизма
$f:{\cal E}\to {\cal F}$ топосов Гротендика имеется коммутативный
квадрат
\[
\xymatrix {
{\bf Geom}({\cal F}, \Set[{\mathbb T}]) \ar[d]^{- \circ f} \ar[rr]^{\simeq}   & & {{\mathbb T}}{\textrm{-mod}}({\cal F})  \ar[d]^{{\mathbb T}\textbf{-mod}(f^{\ast})}  \\
{\bf Geom}({\cal E}, \Set[{\mathbb T}])  \ar[rr]^{\simeq} & &   {\mathbb T}{\textrm{-mod}}({\cal E})}
\]
в (мета-)категории категорий $\textbf{CAT}$.

Напомним, что геометрический морфизм --- естественное,
топологически мотивированное понятие морфизма между топосами
Гротендика; в самом деле, сопоставление, переводящее локаль $L$ в
топос пучков $\Sh(L)$ на ней, поднимается до полного и точного
(faithful) функтора из категории локалей в категорию топосов
Гротендика и геометрических морфизмов между ними, который
идентифицирует первую категорию как полную рефлективную
подкатегорию последней.

Классифицирующий топос геометрической теории $\mathbb T$ может
рассматриваться как \emph{представляющий объект} для
(псевдо-)функтора $\mathbb T$-модели
\[
{{\mathbb T}\textrm{-mod}:\mathfrak{BTop}^{\textrm{op}} \to \textbf{CAT}}
\]
из категории, двойственной к категории $\mathfrak{BTop}$ топосов
Гротендика, в (мета-)категорию $\textbf{CAT}$, который
сопоставляет
\begin{itemize}
\item топосу $\cal E$ категорию ${\mathbb T}\textbf{-mod}({\cal
E})$ моделей для $\mathbb T$ в $\cal E$ и \item геометрическому
морфизму $f:{\cal E}\to {\cal F}$ функтор\\ ${\mathbb
T}\textbf{-mod}(f^{\ast}):{\mathbb T}\textbf{-mod}({\cal F}) \to
{\mathbb T}\textbf{-mod}({\cal E})$,\\ переводящий модель ${M\in
{\mathbb T}\textbf{-mod}({\cal F})}$ в ее образ $f^{\ast}(M)$
относительно функтора обратного образа $f^{\ast}$ для $f$.
\end{itemize}
В частности, классифицирующий топос \emph{единственен с точностью
до эквивалентности категорий}.

Непосредственно связанной с концепцией классифицирующего топоса
является концепция универсальной модели. \emph{Универсальная
модель} геометрической теории $\mathbb T$ есть модель
\emph{$U_{\mathbb T}$} для $\mathbb T$ в топосе Гротендика $\cal
G$ такая, что для любой $\mathbb T$-модели $M$ в топосе Гротендика
$\cal F$ существует единственный (с точностью до изоморфизма)
геометрический морфизм $f_{M}:{\cal F}\to {\cal G}$ такой, что
$f_{M}^{\ast}(U_{\mathbb T})\cong M$.

Согласно ($2$-мерной) лемме Йонеды, если топос $\cal G$ содержит
\emph{универсальную модель} геометрической теории $\mathbb T$, то
$\cal G$ удовлетворяет универсальному свойству
\emph{классифицирующего топоса} для $\mathbb T$. Обратно, если
топос $\cal E$ классифицирует геометрическую теорию $\mathbb T$,
то $\cal E$ содержит универсальную модель для $\mathbb T$.

Следовательно, универсальные модели, как и классифицирующие
топосы, однозначно определены с точностью до категорной
эквивалентности. В частности, если $M$ и $N$ --- универсальные
модели геометрической теории $\mathbb T$, лежащие соответственно в
топосах $\cal F$ и $\cal G$, то существует единственная (с
точностью до изоморфизма) геометрическая эквивалентность между
$\cal F$ и $\cal G$ такая, что ее функторы обратного образа
переводят $M$ и $N$ друг в друга (с точностью до изоморфизма).

Классифицирующие топосы произвольных теорий первого порядка,
вообще говоря, не существуют по той специфической причине, что
функторы обратного образа для геометрических морфизмов могут не
сохранять интерпретацию импликаций и универсальных квантификаций
(даже несмотря на то, что они всегда сохраняют интерпретацию
геометрических формул). С другой стороны, классифицирующие топосы
всегда существуют для геометрических теорий, и они могут быть
построены канонически из таковых по смыслу синтаксической
конструкции (см. раздел \ref{onetopos} ниже). Этот факт имеет
фундаментальное значение, по крайней мере для наших целей, и мы
широко обсудим его ниже.

Заметим, что простое \emph{существование} классифицирующих топосов
является само по себе фактом первичной концептуальной важности;
факт, что функтор $\mathbb T$-модели представим, означает, что по
существу вся информация о функторе или, другими словами, о
категориях моделей теории $\mathbb T$, сосредоточена в
одном-единственном объекте, а именно в классифицирующем топосе.
Это ясно представляет результат `симметрии' по отношению к теории
моделей геометрических теорий, который показывает, что среда
топосов Гротендика предоставляет сильную форму `внутренней
полноты' по отношению к геометрическим теориям: все модели данной
геометрической теории в топосах Гротендика, включая все
классические основанные на множествах модели этой теории, являются
образами единственной `универсальной модели', лежащей в
классифицирующем топосе, в которой имеет место объединение
синтаксических и семантических аспектов теории. Если мы
ограничиваемся рассмотрением моделей теории $\mathbb T$,
основанных на множествах, мы не имеем аналогичного результата
представимости: мы должны расширить наше рассмотрение, например,
до мира топосов Гротендика, чтобы найти такую разновидность
`симметрии'. В конечном счете такой вид явлений является довольно
общим для математики; подумаем, например, что до введения
комплексной плоскости аналитики должны были тратить изрядное время
на попытки понять поведение решений полиномиальных уравнений на
вещественной прямой, тогда как после обнаружения для
полиномиальных уравнений `естественной среды' все начало
осознаваться в новом свете, и люди немедленно перестали
удивляться, почему поведение полиномиальных уравнений на
вещественной прямой так странно. Фактически основная теорема
алгебры по-настоящему воплощает эту достигнутую `симметрию' по
отношению к полиномиальным уравнениям.

Знакомый образ, который приходит на ум при мысли об этом ---
\emph{солнце} и порождаемые им \emph{тени}. Тени возникают, когда
свет солнца встречает некую форму твердой материи; похожим образом
модели возникают, когда фрагмент синтаксиса интерпретируется в
данной `конкретной' среде. Например, абстрактное, синтаксическое
понятие группы порождает много различных моделей в различных
категориях, а именно классическое понятие группы (если мы
интерпретируем ее в категории множеств), понятие топологической
группы (если мы интерпретируем ее в категории топологических
пространств), понятие алгебраической группы (если мы
интерпретируем ее в категории алгебраических многообразий),
понятие группы Ли (если мы интерпретируем ее в категории гладких
многообразий) и т.д.

Как тени легко сравнимы друг с другом, если осознать, что они
происходят от одного источника, так и изучение моделей теории
сильно облегчается рассмотрением синтаксических аспектов теории.
Например, все абстрактные алгебраические свойства групп (в смысле
алгебраических тождеств, доказуемых в аксиоматической теории
групп) могут быть канонически интерпретированы во всех
представленных выше категориях, `автоматически' давая новые
результаты о топологических группах, алгебраических группах,
группах Ли и т.д. Все эти понятия используют общее ядро, которое
находится скорее на синтаксическом уровне, чем на семантическом.

Как солнце является унифицирующим источником теней, так и
синтаксис играет роль унифицирующей концепции для математических
структур. Возможно, это не очень видно в традиционной постановке
классической теории моделей (первого порядка), где имеется теорема
полноты для конечной логики первого порядка, дающая
эквивалентность синтаксической точки зрения и семантической,
базирующейся на теории множеств. Но когда доступно много семантик,
как в теории моделей, основанной на топосах, синтаксис ясно
приобретает особую роль.

\section{Логика, лежащая в основе топосов Гротендика}\label{logGrothendieck}

Мы знаем, что каждая геометрическая теория над заданной сигнатурой
имеет классифицирующий топос, который однозначно определен с
точностью до категориальной эквивалентности. Это естественно
приводит к вопросу, является ли данное сопоставление
разновидностью инъекции или сюръекции, что приводит к
возникновению следующих вопросов:\\

(1) Могут ли две различные геометрические теории иметь эквивалентные классифицирующие топосы?\\

(2) Является ли каждый топос Гротендика классифицирующим топосом некоторой геометрической теории?\\

Ответы на оба эти вопроса хорошо известны. Конкретно, из
характеризации классифицирующего топоса как представляющего
объекта для псевдофунктора моделей (см. раздел \ref{classif})
немедленно следует, что две геометрические теории имеют
эквивалентные классифицирующие топосы, если и только если они
имеют эквивалентные категории моделей в каждом топосе Гротендика
$\cal E$ естественно в $\cal E$. Две такие теории называются
\emph{Морита-эквивалентными}. (Мы вернемся к теме
Морита-эквивалентности в разделе \ref{morita} ниже.) Это ответ на
первый вопрос.

Чтобы ответить на второй вопрос, мы должны проникнуть несколько
глубже в структуру топосов Гротендика. По определению топос
Гротендика есть произвольная категория, эквивалентная категории
$\Sh({\cal C}, J)$ пучков на сайте $({\cal C}, J)$ (сайт
определения топоса Гротендика $\cal E$
--- произвольный сайт $({\cal C}, J)$ такой, что категория
$\Sh({\cal C}, J)$ пучков на $({\cal C}, J)$ эквивалентна $\cal
E$); напомним, что сайт $({\cal C}, J)$ состоит из малой категории
$\cal C$ и топологии Гротендика $J$ на $\cal C$ (мы направляем
читателя к \cite{MM} для первого введения в теорию топосов).

Заметим, что здесь имеется элемент `неканоничности', в котором мы
не можем ожидать, что топос $\Sh({\cal C}, J)$ однозначно
определяет сайт $({\cal C}, J)$; однако этот аспект
`неканоничности' во многих отношениях вообще не является
нежелательным, и фактически он представляет собой фундаментальный
ингредиент взгляда на топосы как на унифицирующие пространства,
описанного ниже.

Напомним, что, для данного сайта $({\cal C}, J)$, для каждого
топоса Гротендика $\cal E$ имеется эквивалентность между
категорией ${\bf Geom}({\cal E}, \Sh({\cal C}, J))$ геометрических
морфизмов из $\cal E$ в $\Sh({\cal C}, J)$ и категорией ${\bf
Flat}_{J}(\cal{C}, \cal{E})$ $J$-непрерывных плоских функторов из
${\cal{C}}$ в ${\cal{E}}$, естественная в $\cal E$. Теперь мы
можем построить геометрическую теорию ${\mathbb T}^{\cal C}_{J}$
такую, что ее модели в любом топосе Гротендика $\cal E$ могут быть
точно отождествлены с $J$-непрерывными плоскими функторами из
${\cal{C}}$ в ${\cal{E}}$ (и гомоморфизмы ${\mathbb T}^{\cal
C}_{J}$-моделей могут быть отождествлены с естественными
преобразованиями между соответствующими плоскими функторами);
понятно, ${\mathbb T}^{\cal C}_{J}$ будет классифицирована топосом
$\Sh({\cal C}, J)$. Мы назовем такую теорию ${\mathbb T}^{\cal
C}_{J}$ \emph{теорией $J$-непрерывных плоских функторов на $\cal
C$}. Это гарантирует, что каждый топос Гротендика возникает как
классифицирующий топос некоторой геометрической теории, давая
положительный ответ на второй вопрос. Полезно выписать явно
аксиоматизацию теории ${\mathbb T}^{\cal C}_{J}$.

Сигнатура для ${\mathbb T}^{\cal C}_{J}$ имеет один и тот же тип
$\name{A}$ для каждого объекта $A$ из $\cal C$, и символ функции
${\name{f}:\name{A} \to \name{B}}$ один для каждой стрелки ${f:A
\to B}$ в $\cal C$. Аксиомы из ${\mathbb T}^{\cal C}_{J}$
следующие (чтобы отметить, что переменная $x$ имеет тип
$\name{A}$, мы пишем $x^{A}$):
\begin{equation}
(\top \vdash_{x} (\name{f}(x)=x))
\end{equation}
для любой тождественной стрелки $f$ в $\cal C$;
\begin{equation}
(\top \vdash_{x} (\name{f}(x)=\name{h}(\name{g}(x))))
\end{equation}
для любой тройки стрелок $f,g,h$ из $\cal C$ такой, что $f$ равна
композиции $h\circ g$;
\begin{equation}
\top \vdash_{[]} \mathbin{\mathop{\textrm{\huge $\vee$}}\limits_{A\in Ob({\cal C})}}(\exists x^{A})\top
\end{equation}
(где дизъюнкция проводится по всем объектам из $\cal C$);
\begin{equation}
(\top \vdash_{x^{A}, y^{B}} \mathbin{\mathop{\textrm{\huge $\vee$}}\limits_{A\stackrel{f}{\leftarrow} C \stackrel{g}{\rightarrow} B}}(\exists z^{C})(\name{f}(z^{C})=x^{A} \wedge \name{g}(z^{C})=y^{B}))
\end{equation}
для любых объектов $A$, $B$ из $\cal C$ (где дизъюнкция проводится
по всем конусам ${A\stackrel{f}{\leftarrow} C
\stackrel{g}{\rightarrow} B}$ на дискретной диаграмме, заданной
парой объектов $A$ и $B$);
\begin{equation}
(\name{f}(x^{A})=\name{g}(x^{A}) \vdash_{x^{A}} \mathbin{\mathop{\textrm{\huge $\vee$}}\limits_{h:C\to A\in Eq(f,g)}}(\exists z^{C})(\name{h}(z^{C})=x^{A}))
\end{equation}
для любой пары стрелок ${f,g:A\to B}$ в $\cal C$ с общими областью
и кообластью (где дизъюнкция проводится по всем стрелкам $h$,
уравнивающим $f$ и $g$);
\begin{equation}
(\top \vdash_{x^{A}} \mathbin{\mathop{\textrm{\huge $\vee$}}\limits_{i\in I}}(\exists y_{i}^{B_{i}})(\name{f_{i}}(y_{i}^{B_{i}})=x^{A}))
\end{equation}
для любого $J$-покрывающего семейства $(f_{i}:B_{i}\to A \textrm{
| } i\in I)$.

Заметим, что первые две группы аксиом выражают функториальность,
третья, четвертая и пятая вместе выражают плоскость (в терминах
фильтрованности соответствующей категории элементов), тогда как
шестая группа аксиом выражает $J$-непрерывность.

Отметим, что даже если $\cal C$ декартова, в каковом случае
плоские функторы на $\cal C$ имеют конечную (когерентную)
аксиоматизацию, наличие топологии Гротендика $J$ делает
аксиоматизацию ${\mathbb T}^{\cal C}_{J}$, вообще говоря,
бесконечной. Это показывает, что `логика, лежащая в основе'
топосов Гротендика, является действительно \emph{геометрической
логикой} во всей своей бесконечности. Дальнейшее свидетельство
этому дает теорема дуальности (см. раздел \ref{duality}), которая
утверждает, что подтопосы классифицирующего топоса геометрической
теории взаимно однозначно соотносятся с \emph{геометрическими}
частными этой теории. Фактически, как мы увидим ниже, большинство
классических теоретико-топосных инвариантов или конструкций ведут
себя вполне естественно по отношению к геометрической логике (это
означает, что они соответствуют естественным свойствам или
операциям на геометрических теориях), но скорее неуклюже по
отношению к \emph{конечному} фрагменту ее, т.е. когерентной
логике. Действительно, многие важные теоретико-топосные
конструкции дают, даже в приложении к когерентным теориям, теории,
которые не являются более когерентными: верный знак, что эти
теоретико-топосно индуцированные преобразования теорий могут
естественно изучаться только в полностью бесконечных рамках.
Другую иллюстрацию существенно геометрической природы логики,
лежащей в основании топосов Гротендика, дает рассмотрение
универсальных моделей геометрических теорий; согласно теореме
7.1.4 \cite{OC}, подобъекты (объекта, лежащего в основе)
универсальной модели геометрической теории могут быть
отождествлены с классами эквивалентности по доказуемости
\emph{геометрических} формул в заданном контексте над сигнатурой
теории.

В этой связи будет дополнительным прояснением напоминание
внутренней характеристики геометрической логики, полученной в
\cite{OC9}:\\

\noindent
 {\bf Теорема~3.1.} {\it Пусть $\Sigma$ --- сигнатура,
$\cal S$
--- набор $\Sigma$-структур в топосах Гротендика, замкнутый
относительно изоморфизмов структур. Тогда $\cal S$ будет набором
всех моделей в топосах Гротендика для геометрической теории над
$\Sigma$, если и только если он удовлетворяет следующим двум
условиям:
\begin{enumerate}[(i)]
\item для любого геометрического морфизма ${f:{\cal F}\to {\cal
E}}$ из принадлежности $M$ набору $\cal S$ следует принадлежность
$f^{\ast}(M)$ набору $\cal S$; \item для любого (индексированного
множеством) совместно сюръективного семейства ${ \{f_{i}:{\cal
E}_{i}\to {\cal E} \textrm{ | } i\in I\} }$ геометрических
морфизмов и любой $\Sigma$-структуры $M$ в $\cal E$ из
принадлежности для каждого ${i\in I}$ $f_{i}^{\ast}(M)$ набору
$\cal S$ следует принадлежность $M$ набору $\cal S$.
\end{enumerate}
}

\section{Один топос --- много сайтов}\label{onetopos}

Классифицирующий топос геометрической теории всегда может быть
построен `канонически' исходя из теории по смыслу синтаксической
конструкции: именно, классифицирующий топос теории $\mathbb T$
дается категорией пучков на \emph{геометрическом синтаксическом
сайте} $({\cal C}_{\mathbb T}, J_{\mathbb T})$ для $\mathbb T$.
Напомним, что \emph{геометрическая синтаксическая категория}
${\cal C}_{\mathbb T}$ геометрической теории $\mathbb T$ над
сигнатурой $\Sigma$ имеет в качестве объектов `классы
эквивалентности по переименованиям' геометрических
формул-в-контексте ${ \{\vec{x}. \phi\} }$ над $\Sigma$ и в
качестве стрелок ${ [\theta]:\{\vec{x}. \phi\} \to \{\vec{y}.
\psi\} }$ классы эквивалентности по $\mathbb T$-доказуемости
геометрических формул $[\theta]$, которые являются $\mathbb
T$-доказуемо функциональными из ${ \{\vec{x}. \phi\} }$ в ${
\{\vec{y}. \psi\} }$, тогда как \emph{геометрическая
синтаксическая топология} $J_{\mathbb T}$ для $\mathbb T$ есть
каноническая топология Гротендика на геометрической категории
${\cal C}_{\mathbb T}$ (см. детали в \cite{El}, часть D).

Если теория $\mathbb T$ лежит в меньшем фрагменте геометрической
логики, таком как декартова, регулярная или когерентная логика,
классифицирующий топос может быть альтернативно вычислен взятием
категории пучков на других синтаксических сайтах, а именно
декартова синтаксическая категория ${\cal C}_{\mathbb
T}^{\textrm{cart}}$ для $\mathbb T$, оснащенная тривиальной
топологией Гротендика (если $\mathbb T$ декартова), регулярный
синтаксический сайт $({\cal C}_{\mathbb T}^{\textrm{reg}},
J_{{\cal C}_{\mathbb T}^{\textrm{reg}}})$ (если $\mathbb T$
регулярна), и когерентный синтаксический сайт $({\cal C}_{\mathbb
T}^{\textrm{coh}}, J_{{\cal C}_{\mathbb T}^{\textrm{coh}}})$ (если
$\mathbb T$ когерентна).

Таким образом, для заданной теории могут быть альтернативные пути
вычисления ее классифицирующего топоса, даже если оставаться в
контексте синтаксических сайтов: например, если теория декартова,
рассмотрев ее как декартову, регулярную, когерентную и
геометрическую теорию, мы получаем четыре различных
`синтаксических сайта' таких, что категория пучков на них дает
классифицирующий топос теории. Этот факт будет использован в
разделе \ref{examp} для порождения различных результатов в логике.

С другой стороны, имеется альтернативный общий метод семантической
природы для вычисления классифицирующих топосов геометрических
теорий, основанный на понятии теории предпучкового типа. Говорят,
что теория имеет предпучковый тип, если она классифицируется
предпучковым топосом. Класс теорий предпучкового типа содержит все
декартовы теории, так же как и многие другие значительные теории
(мы отсылаем читателя к разделу \ref{presheaf} для более широкой
дискуссии на тему теорий предпучкового типа). Классифицирующий
топос теории предпучкового типа $\mathbb T$ дается категорией
функторов $[\textrm{f.p.} {\mathbb T}\textrm{-mod}(\Set), \Set]$,
где $\textrm{f.p.} {\mathbb T}\textrm{-mod}(\Set)$ --- категория
(представлений классов изоморфизмов) конечно представимых моделей
для $\mathbb T$.

Теперь, если геометрическая теория является частным (т.е. теорией,
полученной добавлением геометрических секвенций над той же самой
сигнатурой) для теории предпучкового типа $\mathbb T$, то ее
классифицирующий топос есть подтопос $\Sh(\textrm{f.p.} {\mathbb
T}\textrm{-mod}(\Set)^{\textrm{op}}, J)$ классифицирующего топоса
$[\textrm{f.p.} {\mathbb T}\textrm{-mod}(\Set), \Set]$ для
$\mathbb T$, и топология Гротендика $J$ может быть вычислена
непосредственно переписыванием аксиом теории в форме, использующей
формулы, представляющие (конечно представимую) модель для $\mathbb
T$ (см. детали в \cite{OC}, глава~5).

Если сравнивать с вышеупомянутым методом синтаксических сайтов,
этот последний метод конструирования классифицирующих топосов
является семантическим по духу, даже несмотря на то, что для любой
теории предпучкового типа $\mathbb T$ категория $\textrm{f.p.}
{\mathbb T}\textrm{-mod}(\Set)$ может быть отождествлена с
категорией, обратной полной подкатегории геометрической
синтаксической категории для $\mathbb T$ (см. теорему 10.3.3
\cite{OC} или раздел \ref{presheaf} ниже).

Конечно, имеется много других методов доказательства того, что
конкретный топос Гротендика классифицирует данную геометрическую
теорию; иногда можно в лоб доказать, что модели геометрической
теории могут быть отождествлены (естественно в любом топосе
Гротендика) с плоскими $J$-непрерывными функторами на малой
категории $\cal C$, и в этом случае можно заключить, что
классифицирующий топос теории есть топос $\Sh({\cal C}, J)$ пучков
на сайте $({\cal C}, J)$. Также, работая с топосами, часто можно
доказать эквивалентность топоса, представленного в терминах одного
сайта, топосу, представленному с использованием другого сайта; и,
поскольку понятие классифицирующего топоса, очевидно, инвариантно
относительно категориальной эквивалентности, это может привести к
многим различным `представлениям' классифицирующего топоса для
данной теории.

Другой источник различных сайтов определения для классифицирующего
топоса геометрической теории получается из рассмотрения частных
для геометрических теорий. Согласно теореме дуальности (см. раздел
\ref{duality} ниже), частные геометрической теории $\mathbb T$
могут быть отождествлены с подтопосами классифицирующего топоса
для $\mathbb T$. Далее, понятие подтопоса есть теоретико-топосный
инвариант (т.е. оно зависит только от топоса и не зависит от
конкретных сайтов определений для него), который ведет себя
естественно по отношению к сайтам (действительно, для любого сайта
$({\cal C}, J)$ подтопосы топоса $\Sh({\cal C}, J)$ взаимно
однозначно соответствуют топологиям Гротендика на $\cal C$,
которые содержат $J$); следовательно, любое конкретное
представление классифицирующего топоса геометрической теории дает
соответствующее представление классифицирующего топоса любого его
частного.

Мы можем думать, что каждый сайт определения классифицирующего
топоса геометрической теории представляет частный аспект теории, а
классифицирующий топос воплощает те существенные особенности
теории, которые инвариантны по отношению к частным
(синтаксическим) представлениям теории, индуцирующим
Морита-эквивалентности на семантическом уровне. Мы вернемся к
этому пункту в разделах \ref{onetopos} и \ref{bridge} ниже.

Подчеркнем, что в реальности для заданного топоса Гротендика может
быть очень много различных сайтов определения (логическая
интерпретация этого факта является центральной для наших целей и
будет обсуждаться в разделе \ref{morita}). Зачастую просто взгляд
на одну и ту же теорию с двух разных точек зрения приводит к двум
разным представлениям ее классифицирующего топоса. Например,
когерентная теория полей может быть ясно представлена как в виде
частного теории коммутативных колец с единицей, так и в виде
частного теории регулярных колец фон Неймана, и каждое из этих
представлений дает другое представление классифицирующего топоса
(как категории пучков на категории, обратной к категории конечно
представленных колец, и как категории пучков на категории,
обратной к категории конечно представленных регулярных колец фон
Неймана). Эта техническая `гибкость' теории топосов во включении и
извлечении математической субстанции очевидно нематериального
опыта `видения одного и того же предмета двумя и более способами'
--- один из наиболее поразительных аспектов теории. Как мы
покажем ниже, невероятное количество информации, уместное для
классической математики, `скрыто' внутри топосов и может быть
извлечено с использованием их различных сайтов определений.

\section{Морита-эквивалентности}\label{morita}

Мы видели в разделе \ref{logGrothendieck}, что две геометрические
теории имеют эквивалентные классифицирующие топосы, если и только
если они Морита-эквивалентны. Заметим, что отношение `быть
Морита-эквивалентными' определяет отношение эквивалентности на
совокупности всех геометрических теорий, и топосы Гротендика могут
быть взяты как \emph{канонические представители} результирующих
классов эквивалентности. Теории, которые Морита-эквивалентны друг
другу, являются, в широком смысле, теориями, которые хотя,
возможно, имеют разное лингвистическое (т.е. синтаксическое)
представление, но совместно используют общее `семантическое ядро',
и это ядро в точности воплощено в их общем классифицирующем
топосе.

Морита-эквивалентность --- общее понятие эквивалентности
математических теорий, которое в математике вездесуще (даже
несмотря на то, что в прошлом не было большого интереса к
идентификации Морита-эквивалентностей, возможно, из-за отсутствия
общей теории, приписывающей центральную важность этому понятию и
демонстрирующей его техническую полезность --- одной из целей
настоящей статьи фактически является защита экстремальной важности
исследований в этой области, см. разделы \ref{bridge} и
\ref{work}). Как простой пример теорий, которые
Морита-эквивалентны, можно взять теорию булевых алгебр и теорию
булевых колец.

С точки зрения логика довольно естественно желать рассматривать
как эквивалентные две математические теории, чьи категории
моделей, основанных на множествах, могут быть отождествлены друг с
другом; различие между этим понятием эквивалентности и понятием
Морита-эквивалентности просто в том, что в последнем случае мы
требуем отождествление моделей двух теорий `перенести' на
произвольный топос Гротендика $\cal E$ естественно в $\cal E$. Это
может казаться на первый взгляд очень суровым ограничением; но
фактически мы можем ожидать, что большинство эквивалентностей
первого рода, возникающих в математической практике, расширяемы до
Морита-эквивалентностей. Причина в том, что при установлении
эквивалентности первого рода обычно используются стандартные
теоретико-множественные конструкции, которые не используют закон
исключенного третьего, и поскольку топос Гротендика ведет себя
логически как `обобщенный универсум множеств', в котором можно
представить большинство теоретико-множественных конструкций с
единственным значительным исключением аргументов, требующих закон
исключенного третьего, естественно ожидать возможности `поднять'
такую эквивалентность в произвольный топос Гротендика так, чтобы
глобально получить Морита-эквивалентность.

Что может, пожалуй, быть немного удивительным --- это что теории
непрерывных плоских функторов действительно играют центральную
роль в существенно логической теме Морита-эквивалентности (см.
ниже), также любопытно, что их аксиоматизации весьма эксцентричны
в сравнении со стандартами классической теории моделей (сигнатуры
с, возможно, бесконечным количеством типов (sorts), бесконечные
дизъюнкции и т.д.). В самом деле, $({\cal C}, J)$ есть сайт
определения классифицирующего топоса для геометрической теории
$\mathbb T$, если и только если $\mathbb T$ Морита-эквивалентна
теории $J$-непрерывных плоских функторов на $\cal C$. В частности,
любая геометрическая теория каноническим образом
Морита-эквивалентна теории $J$-непрерывных плоских функторов на
$\cal C$, а именно теории $J_{\mathbb T}$-непрерывных плоских
функторов на ${\cal C}_{\mathbb T}$.

Понятие Морита-эквивалентности также прямо связано с понятием
биинтерпретируемости в классической теории моделей: фактически две
геометрические теории Морита-эквивалентны, если и только если они
биинтерпретируемы одна в другой в обобщенном смысле
(геометрический морфизм между классифицирующими топосами двух
геометрических теорий может рассматриваться как обобщенная
интерпретация одной теории в другой, см. \cite{OC}, раздел 2.1.5).
Если теории когерентны, это в точности то же, что сказать, что их
синтаксические предтопосы эквивалентны, что опять есть обобщенное
понятие биинтерпретируемости.

Тема Морита-эквивалентности прямо соотносится с существованием
различных сайтов определения для данного топоса. Действительно,
два сайта $({\cal C}, J)$ и $({\cal C}', J')$ поднимаются до
одного и того же топоса (т.е. $\Sh({\cal C}, J)$ эквивалентен
$\Sh({\cal C}', J')$), если и только если теории
${J\mbox{-непрерывных}}$ плоских функторов на $\cal C$ и
$J'$-непрерывных плоских функторов на ${\cal C}'$
Морита-эквивалентны; с другой стороны, различные математические
теории, которые Морита-эквивалентны друг другу, дают различные
сайты определений для своих классифицирующих топосов, а именно их
синтаксические сайты. В свете этой связи с понятием сайта и
обсуждений в разделе \ref{onetopos} можно сказать, что понятие
Морита-эквивалентности действительно отлавливает многие
интуитивные идеи `видения одного и того же различными способами'.
(Заметим, что \emph{одна} теория генерирует бесконечное количество
Морита-эквивалентностей: как было отмечено в разделе
\ref{onetopos}, просто взгляд на теорию как на частное той или
другой теории приводит к новому сайту определения для ее
классифицирующего топоса, т.е. к Морита-эквивалентности.)
Фактически данное математическое свойство может проявляться в
нескольких различных формах в контексте математических теорий,
имеющих общее `семантическое ядро', но разные лингвистические
представления; замечателен факт, что если свойство формулируется
как теоретико-топосный инвариант на некотором топосе, то выражение
этого в термах различных теорий, классифицируемых этим топосом,
определяется в значительной степени техническим соотношением между
топосом и различными сайтами определения для него (см. раздел
\ref{bridge}).

В связи с этим представляется полезным отметить, что люди имеют
естественную тенденцию \emph{визуализировать} `семантику' и
использовать `синтаксис' для \emph{лингвистического} обоснования
ее; поэтому `видение одного и того же разными способами' может с
успехом иметь смысл `описания данной структуры с использованием
различных языков'. Фактически в математике (как и в реальной
жизни) часто бывает, что языки и методы, используемые для изучения
некоторого конкретного объекта, могут настолько сильно различаться
от человека к человеку, что довольно трудно \emph{понять}, что
фактически реальный объект изучения один и тот же.

Суммируя: тот факт, что различные теории Морита-эквивалентны одна
другой, переводится теоретико-топосно в существование различных
сайтов определения для классифицирующего топоса. В этом месте
можно естественно задаться вопросом, будет ли эта связь с понятием
сайта сколько-нибудь полезна для исследования
Морита-эквивалентностей. В действительности главная цель этой
статьи --- дать положительный ответ на этот вопрос с акцентом на
наших новых методологиях теоретико-топосной природы для
исследования Морита-эквивалентностей. Эти методологии, которые мы
представляем по ходу статьи, основаны на видении топосов как
`мостов', которые могут быть использованы для переноса информации
между теориями, которые Морита-эквивалентны друг другу, и понятие
сайта играет центральную роль в этом видении (см. раздел
\ref{bridge} ниже). Фактически эти методы дают множество техник
для `унификации математики' в смысле установления новых связей
между разными математическими теориями и \emph{переноса} идей и
результатов между ними.

\section{Топосы как `мосты'}\label{bridge}

Мы уже отмечали, что у данного топоса Гротендика может быть много
различных сайтов определения, и что это соответствует на
логическом уровне существованию Морита-эквивалентностей между
теориями, классифицируемыми данным топосом. Поэтому, в то время
как сопоставление топоса $\Sh({\cal C}, J)$ сайту $({\cal C}, J)$
есть действительно канонический процесс, нахождение теории,
классифицируемой данным топосом Гротендика, не является
каноническим вообще, поскольку это соответствует нахождению малого
сайта определения для топоса: многие различные сайты
сопоставляются посредством конструкции пучка одному и тому же
топосу (с точностью до категорной эквивалентности). Операция
${({\cal C}, J)\rightarrow \Sh({\cal C}, J)}$ взятия пучков на
данном сайте, таким образом, оказывается разновидностью
`кодирования', которое извлекает в точности те существенные детали
теорий, классифицируемых данным топосом, которые инвариантны
относительно Морита-эквивалентности. В известном смысле
классифицирующие топосы воплощают `общие черты' геометрических
теорий, которые эквивалентны друг другу.

Фактически ввиду вышеупомянутых рассмотрений даже не имеет смысла
искать `привилегированный' сайт определения для данного топоса
Гротендика, так как это соответствовало бы способу канонического
выбора теории вне класса Морита-эквивалентности, а это
\emph{заведомо} иррациональное требование, так как не существует
причины, по которой в общем случае нужно было бы предпочесть одну
теорию другой (вообще говоря, очевидно не имеет смысла утверждать
`превосходство' одной ветви математики над другой: все, что можно
рационально сказать --- это что некая разновидность языка может
быть более подходящей, чем другая, \emph{в данном контексте}, но
это субъективные и произвольные соображения).

Конкретно, так как классифицирующий топос геометрической теории
$\mathbb T$ может быть выбран каноническим представителем класса
эквивалентности теорий, которые Морита-эквивалентны $\mathbb T$
(см. раздел \ref{onetopos}), свойства $\mathbb T$, которые
инвариантны по отношению к Морита-эквивалентности, являются, по
крайней мере \emph{концептуально}, свойствами классифицирующего
топоса для $\mathbb T$; обратно, любое свойство классифицирующего
топоса для $\mathbb T$ поднимается до свойства теорий, им
классифицируемых, которые инвариантны относительно
Морита-эквивалентности. \emph{Технически}, исходя из богатства и
гибкости теоретико-топосных методов, мы можем ожидать, что те
свойства геометрических теорий, которые инвариантны относительно
Морита-эквивалентности, в большинстве случаев выражаемы как
инвариантные свойства их классифицирующих топосов, записанные на
теоретико-топосном языке.

Далее, фундаментальной является следующая идея: если мы способны
выразить свойство данной геометрической теории как свойство ее
классифицирующего топоса, то мы можем ожидать выражения этого
свойства в терминах любой другой теории, имеющей тот же самый
классифицирующий топос, и, следовательно, получения соотношения
между исходным свойством и новым свойством другой
Морита-эквивалентной теории. Классифицирующий топос, таким
образом, выступает в роли `моста', связывающего различные
математические теории, Морита-эквивалентные друг другу, что может
использоваться для переноса информации и результатов из одной
теории в другую. Цель настоящей статьи --- показать, что эта идея
топосов как унифицирующих пространств технически вполне
осуществима; великое множество и разнообразие результатов
диссертации \cite{OC} дает ясное свидетельство плодотворности этой
точки зрения, и, как мы показываем по ходу статьи, огромное
количество новых проникновений в любую область математики может
быть получено как результат применения этой техники.

Действительно, факт, что различные математические теории могут
иметь эквивалентные классифицирующие топосы, переходит в
существование различных сайтов определения для одного топоса.
Теоретико-топосные инварианты (т.е. свойства топосов, инвариантные
относительно категорной эквивалентности) могут тогда
использоваться для переноса свойств из одной теории в другую. Как
мы увидим в разделе \ref{invar}, это становится возможным
благодаря тому факту, что абстрактные соотношения между сайтом
$({\cal C},J)$ и топосом $\Sh({\cal C}, J)$, который он
`порождает', часто очень естественны (в том смысле, что свойства
сайтов технически соотносятся с теоретико-топосными инвариантами
естественным образом), давая нам возможность легко передавать
инварианты через различные сайты. Заметим, что, так как построение
категории пучков $\Sh({\cal C}, J)$ по сайту $({\cal C},J)$
полностью канонично, свойство топоса $\Sh({\cal C}, J)$
\emph{является}, по крайней мере в принципе, свойством сайта
$({\cal C},J)$. На практике эти свойства сайтов часто имеют
искренне `категорное' описание, или по меньшей мере
подразумевается или предполагается свойство, допускающее такое
описание. Например, если $\cal C$ --- (малая) категория,
удовлетворяющая правому условию Ора (the right Ore condition), и
$J$ --- атомная топология на ней, то топос $\Sh({\cal C}, J)$
атомный (заметим, что атомность --- теоретико-топосный инвариант).
Для многих инвариантов (например, для свойства топоса быть
булевым, топосом де Моргана или двузначным) имеется биективная
характеризация вида ``$\Sh({\cal C}, J)$ удовлетворяет инварианту,
если и только если сайт $({\cal C}, J)$ удовлетворяет
определенному `приемлемому' категорному свойству'' для
произвольного сайта $({\cal C}, J)$, которая допускает прямой
перенос информации между различными сайтами определения того же
самого топоса. Для других инвариантов могут быть импликации,
проходящие в общем в одном направлении, тогда как для прохождения
в другом направлении должно подразумеваться, что сайт имеет
специальный вид (например, субканонический); как мы увидим в
разделе \ref{logic} ниже, геометрические синтаксические сайты
геометрических теорий ведут себя особенно хорошо по отношению к
этой характеризации. Конечно, может случиться так, что нельзя
установить для данного инварианта общую `приемлемую'
характеризацию сайта вышеупомянутого вида, но в отдельных
интересных случаях можно использовать аргументы \emph{ad hoc} для
установления свойств сайта, из которых следует или которые
являются следствием того факта, что топос удовлетворяет данному
инварианту. Мы увидим эти методологии в действии в многообразии
различных контекстов по ходу статьи.

Как мы увидим в разделах \ref{invar} и \ref{duality} ниже, уровень
общности, представляемый теоретико-топосными инвариантами, идеален
для отлавливания некоторых важных деталей математических теорий.
Дейтвительно, теоретико-топосные инварианты, рассматриваемые на
классифицирующем топосе $\Set[{\mathbb T}]$ геометрической теории
$\mathbb T$, переходят в интересные логические (т.е.
синтаксические или семантические) свойства теории $\mathbb T$, и
теоретико-топосные конструкции на классифицирующих топосах
соответствуют естественным операциям на теориях, классифицируемых
ими.

Суммируя, по смыслу выражения теоретико-топосного инварианта в
терминах различных сайтов определения для данного топоса получаем,
что перенос информации между теориями, классифицируемыми этим
топосом, имеет место; именно в этом смысле топосы действуют как
`мосты', связывающие различные математические теории. Понятие
Морита-эквивалентности, таким образом, приобретает в нашем
контексте выдающуюся роль; в этом смысле оно является примитивным
ингредиентом, на котором только что описанный механизм может быть
запущен. Следовательно, было бы естественным задаться вопросом,
были ли в процессе математической работы обнаружены
Морита-эквивалентности. В действительности имеется много способов
прийти к Морита-эквивалентностям. Мы уже видели, что теория
топосов сама по себе есть первичный источник
Морита-эквивалентностей, так как любой альтернативный способ
представления топоса как категории пучков на сайте приводит к
Морита-эквивалентности (см. разделы \ref{onetopos} и \ref{morita}
выше). Логический подход к Морита-эквивалентностям и задача
расширения классической эквивалентности до Морита-эквивалентности
уже обсуждались в разделе \ref{morita}. Пучковые представления для
различных видов структур должны также допускать естественную
теоретико-топосную интерпретацию в виде Морита-эквивалентностей.
Более того, имеет смысл ожидать, что большинство классических
\emph{дуальностей}, возникающих в математике, использующих некую
геометрическую теорию, должно каким-либо способом подниматься до
Морита-эквивалентностей, извлекающих из них существенные детали.

Естественный способ, которым `работающий математик' может войти в
мир топосов так, чтобы получить пользу из имеющегося знания о
Морита-эквивалентностях и из применения вышеупомянутых методов,
следующий: когда некто находит интересующее его математическое
свойство, ему стоит попытаться найти сайт $({\cal C}, J)$ и
инвариант топоса $\Sh({\cal C}, J)$, который имеет отношение (в
том смысле, что оно подразумевает или подразумевается, хотя бы при
дополнительных допущениях) к исходному свойству. Например, для
данной малой категории ${\cal C}$, удовлетворяющей свойству
амальгамы (см. раздел \ref{fraisse} ниже), свойство категории
${\cal C}$ удовлетворять свойству совместного погружения (см.
раздел \ref{fraisse} ниже), как легко видеть, эквивалентно
свойству топоса $\Sh({\cal C}^{\textrm{op}}, J_{at})$, где
$J_{at}$ --- атомная топология на ${\cal C}^{\textrm{op}}$, быть
двузначным. Так как классифицирующий топос геометрической теории
двузначен, если и только если теория полна (в том смысле, что
любое геометрическое высказывание над сигнатурой теории доказуемо
эквивалентно $\top$ или $\bot$, но не двум им сразу), заключаем,
что свойство совместного погружения на категории $\cal C$, которое
является свойством с `\emph{геометрическим}' ароматом, переходит в
\emph{логическое} свойство полноты теории
$J_{at}\mbox{-непрерывных}$ плоских функторов на ${\cal
C}^{\textrm{op}}$. Этот пример показывает, что общие
математические свойства (в этом смысле --- свойство совместного
погружения на категории $\cal C$) вполне могут возникать как
специализации абстрактных логических свойств теорий, устойчивых
относительно Морита-эквивалентности (в этом смысле --- свойство
полноты геометрических теорий) в виде конкретной геометрической
теории (в этом смысле --- $J_{at}$-непрерывные плоские функторы на
${\cal C}^{\textrm{op}}$). Мы вернемся к этому примеру при
обсуждении теоретико-топосных интерпретаций конструкции Фрессе в
разделе \ref{fraisse}. Отметим центральную роль, которую играют
здесь теории непрерывных плоских функторов; они действительно
имеют сильную связь с классической математикой, и частично причина
заключается в том факте, что их аксиоматизации прямо используют
объекты и стрелки данной категории как соответственно типы (sorts)
и символы функций в своих сигнатурах.

В общем, как можно видеть, придумывание аргументов, подобных
вышеописанному, не настолько трудно и будет становиться со
временем все легче и легче по мере того, как теоретики-топософы
будут открывать новые инварианты к пользе математиков. Как бы то
ни было, мнение автора таково, что это должно было бы быть
совместным усилием; как теоретики-топософы в конструировании
топологических свойств топосов во многом черпали вдохновение из
общей топологии, так в будущем они могли бы быть мотивированы в
своей работе нуждами людей, работающих в любой математической
области. Как мы увидим в разделе \ref{invar} ниже, существование
теоретико-топосных инвариантов с особыми свойствами может иметь
важные ответвления в специфические математические контексты;
топосы действительно играют центральную роль в математике, и их
использование может быть стратегическим в тяжелых ситуациях.

Использование только что описанных методологий позволяет извлечь
невероятное количество новой информации о Морита-эквивалентностях
и установить такие связи между различными теориями, которые с
большим трудом могут быть обнаружены другими способами. Фактически
разновидность понимания, которую эти методы могут порождать,
существенно отличается от (если не принадлежит на самом деле к;
пожалуй, говорить об этом рановато) той, что дается традиционными
методами переноса информации между эквивалентными теориями с
использованием специфического описания эквивалентности.
Действительно, в нашем подходе вместо использования явного
описания Морита-эквиваентности разрабатывается соотношение между
топосом и его сайтами определения. Это возможно, потому что у нас
есть актуальный математический объект, а именно классифицирующий
топос, который извлекает `общее ядро' теорий, им классифицируемых,
и эквивалентности теорий существенным образом закодированы в
соотношении между топосом и его различными сайтами определения.
Фактически для большинства целей реально имеет значение только
\emph{существование} Морита-эквивалентности, и мы с успехом можем
игнорировать актуальное описание ее; конечно, если хочется
установить более `специфические' результаты, явное описание
Морита-эквивалентности становится необходимым (в этом смысле
скорее можно использовать инвариантные свойства \emph{объектов} из
топосов, чем инвариантные свойства `целых топосов', см. раздел
\ref{examp} ниже), но для создания большинства `глобальных'
свойств теорий это не является необходимым вообще (так как по
определению теоретико-топосный инвариант устойчив относительно
\emph{любого} рода категориальной эквивалентности). Другими
словами, мы можем генерировать большое количество интересных
прозрений, просто оставаясь на `одномерном уровне'. Подчеркнем,
что в только что разъясненной технике имеется значительный элемент
\emph{автоматизма}; по смыслу этих методов можно генерировать
новые математические результаты без реальных усилий осознания:
действительно, в большинстве случаев можно просто с готовностью
применить хорошо известные характеризации, связывающие свойства
сайтов и теоретико-топосных инвариантов (таких, как, например,
вышеупомянутые характеризации для свойств булевости и
двузначности) к конкретному интересующему случаю. С другой
стороны, уровень применимости этих методов в математике
неограничен ввиду большой общности понятия топоса.

Исследование теорий предпучкового типа доведено до конца в
диссертации автора и кратко обсуждается в разделе \ref{presheaf}
ниже, представляющем ясную иллюстрацию этих принципов. Фактически
если теория имеет предпучковый тип, мы автоматически имеем два
различных представления ее классифицирующего топоса:
$[\textrm{f.p.} {\mathbb T}\textrm{-mod}(\Set), \Set]$ и
$\Sh({\cal C}_{\mathbb T}, J_{\mathbb T})$. Исходя из данного
двойного представления классифицирующего топоса, можно получить
огромное число глубоких результатов, используя эти методологии,
включая общий вариант теоремы Фрессе в классической теории моделей
(см. раздел \ref{fraisse}).

\section{Теоретико-топосные инварианты}\label{invar}

Мы уже объяснили роль теоретико-топосных инвариантов в нашем
взгляде на топосы как на унифицирующие пространства, которые могут
действовать как `мосты', связывающие различные математические
теории, Морита-эквивалентные друг другу. Но что мы в точности
подразумеваем под выражением `теоретико-топосный инвариант'? Под
этим выражением мы понимаем любое \emph{свойство} $P$ для
(семейств) топосов Гротендика или \emph{конструкцию} $C$,
использующую топосы Гротендика (или их семейства), которые
инвариантны относительно категорной эквивалентности топосов
Гротендика, т.е. такие, что если данное семейство топосов
Гротендика ${\{{\cal E}_{i} \textrm{ | } i \in I\}}$ удовлетворяет
свойству $P$, то для любого такого семейства ${\{{\cal E}_{i}'
\textrm{ | } i \in I\}}$ топосов, что для каждого ${i\in I}$ топос
${\cal E}_{i}$ эквивалентен топосу ${\cal E}_{i}'$, последнее
семейство также удовлетворяет свойству $P$ (требование к
конструкции $C$ таково, что если ${\{{\cal E}_{i} \textrm{ | } i
\in I\}}$ и ${\{{\cal E}_{i}' \textrm{ | } i \in I\}}$ --- два
семейства топосов такие, что для каждого ${i\in I}$ топос ${\cal
E}_{i}$ эквивалентен топосу ${\cal E}_{i}'$, то результат
применения конструкции $C$ к первому семейству должно быть
эквивалентно --- в математическом смысле --- результату применения
$C$ ко второму семейству). Мы \emph{не} требуем, чтобы имелось
лингвистическое выражение свойства $P$ (соответственно описание
конструкции $C$) на языке теории топосов; все, что имеет значение
--- это что $P$ должно быть инвариантно относительно категорной
эквивалентности топосов. Конечно, любое свойство, которое
допускает лингвистическое описание на (неформальном) языке теории
топосов, автоматически является теоретико-топосным инвариантом
(как следствие очень общей `мета-теоремы'), но не следует
ограничиваться этим соображением в поиске инвариантов. Фактически
мы можем ожидать вокруг гораздо большего числа теоретико-топосных
инвариантов, чем те, которые мы можем в настоящее время описать на
категорном языке, обычно используемом для разговоров и рассуждений
о топосах.

Примеры хорошо известных теоретико-топосных инвариантов включают:
свойство топоса быть \emph{булевым}, топосом \emph{де~Моргана},
\emph{атомным}, \emph{двузначным}, быть \emph{связным}, быть
\emph{локально связным}, быть \emph{компактным}, быть
\emph{локальным}, быть \emph{подтопосом} данного топоса (в смысле
геометрического включения), быть \emph{классифицирующим топосом}
данной геометрической теории, иметь достаточно точек и т.д. Группы
когомологий и гомотопий топосов также являются важными
теоретико-топосными инвариантами.

Как мы предполагали в предыдущем разделе, большинство этих
инвариантов ведет себя вполне естественно по отношению к сайтам, а
это означает, что у нас есть общая характеризация, связывающая
`естественные' свойства сайтов с этими инвариантами на
соответствующих топосах пучков. Фактически для большинства
`топологически мотивированных' инвариантов топосов у нас есть
естественные характеризации, идущие в одном направлении, т.е.
утверждения импликаций вида ``если сайт $({\cal C}, J)$
удовлетворяет данному свойству, тогда топос $\Sh({\cal C}, J)$
удовлетворяет данному инварианту'' (мы сошлемся на \cite{El},
часть C, детально представляющую эти результаты), и при
специфических условиях на сайт часто можно установить результаты,
идущие в обратном направлении (см. раздел \ref{logic} ниже). С
другой стороны, для других `логически мотивированных' инвариантов,
таких как свойство топоса быть булевым, топосом де Моргана или
двузначным, мы имеем естественные \emph{биективные} характеризации
(см. \cite{OC}, в особенности главы $6$ и $9$). В разделе
\ref{duality} мы сфокусируемся на конкретном инварианте, а именно
на понятии \emph{подтопоса}, и обсудим его отношение к
классической математике.

Как мы уже отмечали в предыдущем разделе, обычно более принято
рассматривать инварианты объектов топосов с целью установления
`локальных' свойств, чем поиска `глобальных'. Под `инвариантом
объектов топосов' мы понимаем свойство $Q$ объектов топосов такое,
что если ${\tau:{\cal E}\to {\cal F}}$ --- эквивалентность топосов
и $a$ --- объект из ${\cal E}$, то объект $a$ удовлетворяет $Q$,
если и только если объект $\tau(a)$ удовлетворяет $Q$ (и
аналогично для семейств объектов). Примеры инвариантов объекта
топосов включают: свойство объекта быть компактным, связным,
неразложимым, неприводимым, когерентным, быть атомом и т.д.
(формальное определение этих инвариантов см. в \cite{OC}, глава
10).

Мы видели, что любой теоретико-топосный инвариант может с помощью
описанных выше методологий генерировать заметное количество
математических результатов. С другой стороны, также и
`отрицательные результаты', утверждающие \emph{несуществование}
теоретико-топосных инвариантов с особыми свойствами, могут быть
прямо привязаны к соответствующим результатам в математике.
Например, можно задаться вопросом, существует ли
теоретико-топосный инвариант $P$ со следующим свойством: если
$\cal E$ --- топос пучков $\Sh_{G}(X)$ на топологическом группоиде
$s,t:G\to X$, то $\cal E$ удовлетворяет $P$, если и только если
группоид $s,t:G\to X$ `алгебраически связен', т.е. коуравнитель
для $s$ и $t$ в категории топологических пространств дается
одноточечным пространством. Инвариант, удовлетворяющий аналогу
этого свойства для локальных группоидов, существует, и фактически
это в точности свойство топоса быть гиперсвязным (см. \cite{El},
лемма C5.3.7). С другой стороны, такой инвариант не может
существовать для топологических группоидов, и несуществование
такого инварианта прямо соотносится с классическим
теоретико-модельным фактом, что могут существовать (когерентные)
теории, которые \emph{не} вычислимо категориальны, но которые
$k$-категориальны для некоторого невычислимого кардинала $k$ (или
наоборот). Действительно, из теоремы представления (I.~Moerdijk,
C.~Butz) топосов Гротендика с достаточным количеством точек как
топосов пучков на топологическом группоиде (см. \cite{BM})
следует, что классифицирующий топос когерентной теории $\cal T$
может быть представлен как топос пучков $\Sh_{G}(X)$ на
топологическом группоиде $s,t:G\to X$ изоморфизмов моделей для
$\mathbb T$ в $X$, где $X$ --- \emph{произвольное} множество
(занумерованных) моделей для $\mathbb T$ в $\Set$, которое
совместно консервативно для $\mathbb T$. Далее, так как по теореме
Левенгейма-Скулема набор всех моделей вычислимой теории
фиксированной мощности (бесконечно вычислимой или невычислимой)
совместно консервативен для этой теории, то существование нашего
инварианта $P$ подразумевает, что любая вычислимая когерентная
теория $\mathbb T$ имеет в точности один класс изоморфизмов
моделей мощности $k$, если и только если она имеет в точности один
класс изоморфизмов моделей мощности $k'$ для \emph{любых}
кардиналов ${k,k' \geq \omega}$, и, как хорошо известно, это
неверно в общем случае (контрпример дает когерентная теория
алгебраически замкнутых полей заданной характеристики).

\subsection{Логический смысл инвариантов}\label{logic}

Ввиду наших предыдущих утверждений, касающихся важности
теоретико-топосных инвариантов для исследования геометрических
теорий, естественно задаться вопросом, соотносятся ли важные
инварианты классифицирующих топосов с интересными логическими
(т.е. синтаксическими или семантическими) свойствами теорий, ими
классифицируемых. Цель предыдущего раздела --- показать, что
соотносятся. Фактически эта линия исследования была систематически
доведена до конца в \cite{OC}, в основном в главах 6 и 10, где
многие хорошо известные инварианты были изучены с точки зрения
геометрических теорий. Основным инструментом, используемым в этих
исследованиях, является геометрический синтаксический сайт
геометрической теории. Напомним, что классифицирующий топос
геометрической теории $\mathbb T$ всегда может быть представлен
как категория $\Sh({\cal C}_{\mathbb T}, J_{\mathbb T})$ пучков на
геометрическом синтаксическом сайте $({\cal C}_{\mathbb T},
J_{\mathbb T})$ для $\mathbb T$. Далее, этот сайт ведет себя
вполне хорошо в отношении задачи переноса свойств от топоса к
сайту. Действительно, топология Гротендика $J_{\mathbb T}$
субканоническая, и вложение Йонеды ${y:{\cal C}_{\mathbb T} \to
\Sh({\cal C}_{\mathbb T}, J_{\mathbb T})}$ удовлетворяет
дополнительному полезному свойству, что любой подобъект в
$\Sh({\cal C}_{\mathbb T}, J_{\mathbb T})$ объекта вида $y(c)$
имеет вид ${y(d)\mono y(c)}$ для некоторого подобъекта ${d\mono
c}$ в $\Sh({\cal C}_{\mathbb T}, J_{\mathbb T})$. Таким образом,
поскольку свойства синтаксического сайта теории часто
перефразируются естественным образом как синтаксические свойства
теории, процесс выражения теоретико-топосных инвариантов на
классифицирующем топосе теории $\mathbb T$ как синтаксических
свойств $\mathbb T$ работает в целом достаточно гладко. Чтобы это
проиллюстрировать, мы изложим ниже некоторые результаты из
\cite{OC}, касающиеся логической интерпретации инвариантов.
Подчеркнем, что эти результаты имеют место \emph{единообразно} для
любой геометрической теории, т.е. для получения такой
характеризации не подразумевается никаких особых свойств теории.
Чтобы представить результаты, мы сначала должны ввести некоторую
терминологию.\\

\noindent
 {\bf Определение 7.1.}
 Пусть $\mathbb T$ --- геометрическая теория над сигнатурой
$\Sigma$, и $\phi(\vec{x})$ --- геометрическая формула-в-контексте
над $\Sigma$. Тогда

\begin{enumerate}[(i)]

\item Мы говорим, что $\phi(\vec{x})$ \emph{$\mathbb T$-полна},
если секвенция ($\phi \vdash_{\vec{x}} \bot$) не является
доказуемой в $\mathbb T$, и для каждой геометрической формулы
$\chi(\vec{x})$ в том же самом контексте либо ($\phi
\vdash_{\vec{x}} \chi$), либо ($\chi \wedge \phi \vdash_{\vec{x}}
\bot$) доказуемо в $\mathbb T$.

\item Мы говорим, что $\phi(\vec{x})$ \emph{$\mathbb
T$-неразложима}, если для любого семейства $\{\psi_{i}(\vec{x})
\textrm{ | } i\in I\}$ геометрических формул в том же самом
контексте таких, что для каждого $i$ $\psi_{i}$ $\mathbb
T$-доказуемо следует из $\phi$ и для любых различных ${i, j\in I}$
$\psi_{i}\wedge \psi_{j} \vdash_{\vec{x}} \bot$ доказуемо в
$\mathbb T$, мы имеем, что
 $$
   \phi \vdash_{\vec{x}}
   \mathbin{\mathop{\textrm{\huge $\vee$}}\limits_{i\in I}} \psi_{i}
   \mbox{ доказуемо в } \mathbb{T}
 $$
следует из
 $$
    \phi \vdash_{\vec{x}} \psi_{i} \mbox{ доказуемо в } \mathbb{T}
    \mbox{ для некоторого } i\in I.
 $$

\item Мы говорим, что $\phi(\vec{x})$ \emph{$\mathbb
T$-несводима}, если для любого семейства ${\{\theta_{i} \textrm{ |
} {i\in I} \}}$ $\mathbb T$-доказуемо функциональных
геометрических формул ${\{\vec{x_{i}}, \vec{x}.\theta_{i}\}}$ из
${\{\vec{x_{i}}. \phi_{i}\}}$ в ${\{\vec{x}. \phi\}}$ таких, что
${\phi \vdash_{\vec{x}} \mathbin{\mathop{\textrm{\huge
$\vee$}}\limits_{i\in I}}(\exists \vec{x_{i}})\theta_{i}}$
доказуема в $\mathbb T$, существуют ${i\in I}$ и $\mathbb
T$-доказуемо функциональная геометрическая формула ${\{\vec{x},
\vec{x_{i}}. \theta'\}}$ из ${\{\vec{x}. \phi\}}$ в
${\{\vec{x_{i}}. \phi_{i}\}}$ такая, что ${\phi \vdash_{\vec{x}}
(\exists \vec{x_{i}})(\theta' \wedge \theta_{i})}$ доказуема в
$\mathbb T$.

\item Мы говорим, что $\phi(\vec{x})$ \emph{$\mathbb
T$-компактна}, если для любого семейства ${\{\psi_{i}(\vec{x})
\textrm{ | } i\in I\}}$ геометрических формул в том же контексте
из
 $$
   \phi \vdash_{\vec{x}} \mathbin{\mathop{\mbox{\huge
   $\vee$}}\limits_{i\in I}} \psi_{i}
   \mbox{ доказуема в } \mathbb{T}
 $$
следует
 $$
   \phi \vdash_{\vec{x}} \mathbin{\mathop{\mbox{\huge
   $\vee$}}\limits_{i\in I'}} \psi_{i}
   \mbox{ доказуема в } \mathbb{T}
 $$
для некоторого конечного подмножества $I'$ из $I$.

\end{enumerate}


Эти понятия возникают в точности из переформулировки инвариантных
свойств объектов из топосов --- а именно свойств объекта быть
атомом, быть неразложимым, быть несводимым, быть компактным --- в
терминах объектов синтаксического сайта теории; как таковые они
представляют собой реализацию абстрактных методологий,
представленных в разделе \ref{bridge}, для случая
классифицирующего топоса геометрической теории, представленного в
терминах синтаксического сайта этой теории. Конкретно, если
${y:{\cal C}_{\mathbb T} \to \Sh({\cal C}_{\mathbb T}, J_{\mathbb
T})}$ --- вложение Йонеды, то геометрическая формула-в-контексте
${\{\vec{x}. \phi\}}$ $\mathbb T$-полна (соответственно $\mathbb
T$-неразложима, $\mathbb T$-несводима, $\mathbb T$-компактна),
если и только если объект ${y(\{\vec{x}. \phi\})}$ в топосе
$\Sh({\cal C}_{\mathbb T}, J_{\mathbb T})$ является атомом
(соответственно неразложимым, несводимым, компактным). Последний
шаг для получения связи с инвариантами классифицирующего топоса
для $\mathbb T$ --- соотнести эти `глобальные' свойства топосов с
вышеупомянутыми `локальными' свойствами объектов; это оказывается
возможным в данных случаях, и результат следующий.\\

\noindent
 {\bf Теорема 7.2.}
{\it
 Пусть $\mathbb T$ --- геометрическая теория над сигнатурой
$\Sigma$, и $\Set[{\mathbb T}]$ --- ее классифицирующий топос.
Тогда

\begin{enumerate}

\item $\Set[{\mathbb T}]$ \emph{локально связен} (соответственно
\emph{атомный}), если и только если для любой геометрической
формулы $\phi(\vec{x})$ над $\Sigma$ существует (единственное)
семейство ${\{\psi_{i}(\vec{x}) \textrm{ | } i\in I\}}$ $\mathbb
T\mbox{-неразложимых}$ (соответственно $\mathbb T$-полных)
геометрических формул в том же самом контексте таких, что

\begin{enumerate}[(i)]
\item для каждого $i$ из $\mathbb T$-доказуемости $\psi_{i}$
следует $\phi$,

\item для любых различных ${i, j\in I}$ ${\psi_{i}\wedge \psi_{j}
\vdash_{\vec{x}} \bot}$ доказуемы в $\mathbb T$, и

\item $\phi \vdash_{\vec{x}} \mathbin{\mathop{\textrm{\huge
$\vee$}}\limits_{i\in I}} \psi_{i}$ доказуема в $\mathbb T$.
\end{enumerate}

\item $\Set[{\mathbb T}]$ \emph{эквивалентен предпучковому
топосу}, если и только если существует набор $\cal F$ из $\mathbb
T$-несводимых геометрических формул-в-контексте над $\Sigma$,
удовлетворяющий следующему свойству: для любой геометрической
формулы ${\{\vec{y}. \psi\}}$ над $\Sigma$ существуют объекты
${\{\vec{x_{i}}. \phi_{i}\}}$ в $\cal F$, $i$ пробегает $I$, и
$\mathbb T$-доказуемо функциональные геометрические формулы
${\{\vec{x_{i}}, \vec{y}.\theta_{i}\}}$ из ${\{\vec{x_{i}}.
\phi_{i}\}}$ в ${\{\vec{y}. \psi\}}$ такие, что ${\psi
\vdash_{\vec{y}} \mathbin{\mathop{\textrm{\huge
$\vee$}}\limits_{i\in I}}(\exists \vec{x_{i}})\theta_{i}}$
доказуема в $\mathbb T$;

\item $\Set[{\mathbb T}]$ \emph{компактен}, если и только если
формула $\top$ в пустом контексте $\mathbb T$-компактна.

\item $\Set[{\mathbb T}]$ \emph{двузначен}, если и только если
формула $\top$ в пустом контексте $\mathbb T$-полна.

\end{enumerate}
}

Другие результаты, касающиеся логической интерпретации
теоретико-топосных инвариантов, которые были получены в \cite{OC},
следующие:\\

\noindent
 {\bf Теорема 7.3.}
{\it
 Пусть $\mathbb T$ --- геометрическая теория над сигнатурой
$\Sigma$, и пусть $\Set[{\mathbb T}]$ --- ее классифицирующий
топос. Тогда

\begin{enumerate}

\item $\Set[{\mathbb T}]$ \emph{булев}, если и только если для
любой геометрической формулы $\phi(\vec{x})$ над $\Sigma$ имеется
геометрическая формула $\chi(\vec{x})$ над $\Sigma$ в том же самом
контексте такая, что ${\phi(\vec{x}) \wedge \chi(\vec{x})
\:\vdash_{\vec{x}}\: \bot}$ и ${\top \:\vdash_{\vec{x}}\:
\phi(\vec{x}) \vee \chi(\vec{x})}$ доказуемы в $\mathbb T$.

\item $\Set[{\mathbb T}]$ является топосом \emph{де Моргана}, если
и только если для любой геометрической формулы $\phi(\vec{x})$ над
$\Sigma$ существуют две согласованные геометрические формулы
$\psi_{1}(\vec{x})$ и $\psi_{2}(\vec{x})$ над $\Sigma$
в том же самом контексте такие, что:\\
$\top \vdash_{\vec{x}} \psi_{1}(\vec{x})\vee \psi_{2}(\vec{x})$ доказуема в $\mathbb T$,\\
$\psi_{1}(\vec{x})\wedge \phi(\vec{x}) \vdash_{\vec{x}} \bot$ доказуема в $\mathbb T$, и\\
для любой геометрической формулы $\chi(\vec{x})$ над $\Sigma$ в
том же самом контексте такой, что ${\chi(\vec{x}) \vdash_{\vec{x}}
\psi_{2}(\vec{x})}$ доказуема в $\mathbb T$, ${\chi(\vec{x})
\wedge \phi(\vec{x}) \vdash_{\vec{x}} \bot}$ доказуема в $\mathbb
T$, если и только если ${\chi(\vec{x}) \vdash_{\vec{x}} \bot}$
доказуема в $\mathbb T$.

\item $\Set[{\mathbb T}]$ имеет \emph{достаточно точек}, если и
только если, когда геометрическая секвенция $\sigma$ над $\Sigma$
действительна в любой ${\Set\mbox{-модели}}$ для $\mathbb T$,
$\sigma$ доказуема в $\mathbb T$.

\end{enumerate}
}

Эти результаты показывают, что теоретико-топосные инварианты
классифицирующего топоса геометрической теории действительно
соответствуют естественным и интересным логическим свойствам
теории `единообразно' для любой геометрической теории; мы увидим
некоторые приложения этих понятий в следующих разделах статьи.

В связи с этим мы также заметим, что, так как каждый топос
Гротендика имеет синтаксический сайт определения (являясь
классифицирующим топосом геометрической теории), у всех
инвариантных свойств или конструкций на топосах, и, в частности, у
всех результатов теории топосов есть, по крайней мере в принципе,
логический двойник; и так как перенос свойств из топоса в
синтаксический сайт и затем в геометрическую теорию работает, в
общем, достаточно гладко, мы можем ожидать для большинства
инвариантов топосов Гротендика поднятия до свойств геометрических
теорий, выраженных на языке геометрической логики. Например, мы
имеем в виду теорему Делиня о когерентных топосах, которая, как
хорошо известно, эквивалентна классической теореме полноты для
когерентной логики (фактически, как мы увидим в разделе
\ref{examp}, эта эквивалентность подпадает под нашу общую схему
для переноса знаний между различными математическими теориями с
использованием топосов).

Итак, (геометрическая) логика дает язык, на котором б\'{о}льшая
часть теории топосов (и, значит, математики) может быть прочтена и
--- мы подчеркиваем --- не только \emph{сформулирована}, но и
действительно \emph{сделана}. Это показывает, что логика в целом и
геометрическая логика в частности могут быть фундаментальными
инструментами для решения огромного многообразия математических
проблем. В своей работе специалисты по теории моделей уже
продемонстрировали фундаментальное воздействие, которое логические
исследования могут оказать на классическую математику; мы верим,
что эта позитивная тенденция будет продолжаться в будущем и будет
значительно обогащена контактом с теоретико-топосной методологией.

\section{Теорема дуальности}\label{duality}

В этом разделе мы обсуждаем важность понятия \emph{подтопоса} для
геометрической логики и математики в целом; как мы увидим ниже,
это теоретико-топосный инвариант, который ведет себя особенно
хорошо по отношению к сайтам; действительно, подтопосы топоса
$\Sh({\cal C}, J)$ биективно соотносятся с топологиями Гротендика
на $\cal C$, которые содержат $J$. Это позволяет нам плодотворно
применять нашу философию `топосы как мосты' в контексте
подтопосов.

В \cite{OC} доказана \emph{теорема дуальности}, связывающая
подтопосы классифицирующего топоса геометрической теории $\mathbb
T$ и геометрические `частные' для $\mathbb T$.

Прежде чем сформулировать теорему, нам нужно ввести пару
определений.\\

\noindent
 {\bf Определение~8.1.}
 Пусть $\mathbb{T}$ --- геометрическая теория над сигнатурой
$\Sigma$. \emph{Частное} для $\mathbb{T}$ --- это геометрическая
теория $\mathbb{T}'$ над $\Sigma$ такая, что каждая аксиома для
$\mathbb{T}$ доказуема в $\mathbb{T}'$.\\

\noindent {\bf Определение~8.2.}
  Пусть $\mathbb{T}$ и $\mathbb{T}'$ --- геометрические теории над
сигнатурой~$\Sigma$. Мы говорим, что $\mathbb{T}$ и $\mathbb{T}'$
\emph{синтаксически эквивалентны}, и пишем ${\mathbb{T} \equiv_{s}
\mathbb{T}'}$, если для любой геометрической секвенции $\sigma$
над $\Sigma$ имеет место утверждение: $\sigma$ доказуема в
$\mathbb{T}$, если и только если $\sigma$ доказуема в
$\mathbb{T}'$.\\

\noindent
 {\bf Теорема~8.3.}
 {\it
 Пусть $\mathbb{T}$ --- геометрическая теория над сигнатурой
$\Sigma$. Тогда соответствие, переводящее частное для $\mathbb{T}$
в его классифицирующий топос, определяет биекцию между классами
${{\equiv_{s}}\mbox{-эквивалентности}}$ частных для $\mathbb T$ и
подтопосами классифицирующего топоса $\Set[\mathbb{T}]$ для
$\mathbb{T}$.
 }\\

Биекция, задаваемая теоремой, естественна в следующем смысле. Если
${i_{{\cal F}}:{\cal F}\hookrightarrow \Set[\mathbb{T}]}$ ---
подтопос для $\Set[\mathbb{T}]$, отвечающий частному ${\mathbb
T}'$ для $\mathbb T$ посредством теоремы дуальности, имеется
коммутативная (с точностью до естественного изоморфизма) диаграмма
в $\textbf{CAT}$ (где $i$ --- очевидное включение)

\[
\xymatrix {
{{\mathbb T}'}\textrm{-mod}({\cal E}) \ar[rr]^{\simeq} \ar[d]^{i} & & {\bf Geom}({\cal E}, {\cal F}) \ar[d]^{i_{{\cal F}}\circ -} \\
{\mathbb T}\textrm{-mod}({\cal E}) \ar[rr]^{\simeq} & &  {\bf Geom}({\cal E},\Set[\mathbb{T}])}
\]

 \noindent
 естественно в ${\cal E}\in \mathfrak{BTop}$.

Заметим, что немедленным следствием теоремы дуальности является
тот факт, что две геометрические теории над одной и той же
сигнатурой имеют эквивалентные классифицирующие топосы, если и
только если они синтаксически эквивалентны.

В \cite{OC} даны два различных доказательства этой теоремы: одно
опирается на теорию классифицирующих топосов, другое основано на
интерпретации понятия топологии Гротендика с помощью теории
доказательств (что мы обсудим ниже). Как бы то ни было, оба
аргумента существенно используют синтаксическое представление
классифицирующего топоса для $\mathbb T$ как категории пучков на
геометрическом синтаксическом сайте для $\mathbb T$. Таким
образом, эти результаты подпадают под обсуждавшуюся в \ref{logic}
общую схему переформулирования теоретико-топосных инвариантов
классифицирующего топоса теории $\mathbb T$ в терминах логических
свойств или конструкций на $\mathbb T$ с помощью пропускания через
геометрический синтаксический сайт для $\mathbb T$.

Теорема дуальности реализует унификацию теории элементарных
топосов с геометрической логикой, пропуская ее через теорию
топосов Гротендика. Действительно, теорема позволяет нам
интерпретировать многие концепции теории элементарных топосов,
которая применяется к решетке подтопосов данного топоса на уровне
геометрических теорий. Эти понятия включают, например, структуру
когейтинговой алгебры на решетке подтопосов данного топоса,
открытые, замкнутые, квазизамкнутые подтопосы, факторизацию
плотность-замкнутость (dense-closed factorization) геометрического
вложения, когерентные подтопосы, подтопосы с достаточным
количеством точек, факторизацию сюръекция-вложение
(surjection-inclusion factorization) геометрического морфизма,
скелетные вложения, атомы в решетке подтопосов данного топоса,
\emph{булеанизацию} и \emph{деморганизацию} топоса. Замечательный
факт: результирующие понятия и результаты в геометрической логике
представляют заметный логический интерес. Ниже мы рассмотрим
некоторые из них.

Во-первых, мы отметим, что теорема дуальности имеет элегантную
интерпретацию в теории доказательств. Чтобы описать ее, мы сначала
должны заметить, что понятие топологии Гротендика на $\cal C$
естественно поднимается до абстрактной `системы доказательств'.
Конкретно, для данного набора $\cal A$ решет на данной категории
$\cal C$ мы можем определить систему доказательств ${\cal T}_{\cal
C}^{\cal A}$ следующим образом: аксиомы для ${\cal T}_{\cal
C}^{\cal A}$
--- это решета в $\cal A$ вместе со всеми максимальными решетами,
тогда как правила вывода для ${\cal T}_{\cal C}^{\cal A}$ --- это
теоретико-доказательные версии хорошо известых аксиом для
топологий Гротендика, т.е. правила:

\emph{Правило устойчивости:}
\[
\begin{array}{c}
R\\
\hline f^{\ast}(R)
\end{array}
\]
Здесь $R$ --- произвольное решето на объекте $c$ в ${\cal C}$, $f$
--- произвольная стрелка в ${\cal C}$ с кообластью $c$.

\emph{Правило транзитивности:}
\[
\begin{array}{c}
Z \textrm{     } \{f^{\ast}(R) \textrm{ | } f\in Z \}\\
\hline R
\end{array}
\]
Здесь $R$ и $Z$ --- решета в ${\cal C}$ на данном объекте из $\cal
C$.

В этих терминах теорема дуальности дает для любой геометрической
теории $\mathbb T$ некую разновидность
\emph{теоретико-доказательной эквивалентности} между традиционной
системой доказательств геометрической логики с добавлением аксиом
из $\mathbb T$ и системой ${\cal T}_{\cal C}^{\cal A}$, где ${\cal
C}$ --- геометрическая синтаксическая категория ${\cal C}_{\mathbb
T}$ теории $\mathbb T$, ${\cal A}$ --- набор решет в
синтаксической топологии $J_{\mathbb T}$ для $\mathbb T$. В
частности, с каждой геометрической секвенцией над сигнатурой для
$\mathbb T$ можно ассоциировать решето в категории ${\cal C}$ ---
и наоборот --- таким способом, что эти два сопоставления
естественны по отношению к понятиям доказуемости двух систем
доказательств и обратны друг другу с точностью до доказуемости.
`Замкнутые' теории первой системы доказательств могут быть ясно
отождествлены с классами $\equiv_{s}$-эквивалентности частных для
$\mathbb T$, тогда как `замкнутые' теории второй системы
доказательств являются топологиями Гротендика на ${\cal
C}_{\mathbb T}$, которые содержат $J_{\mathbb T}$ (за
дополнительными деталями мы направляем читателя к \cite{OC},
глава~2).

Похожие эквивалентности между классическими системами
доказательств геометрической логики и системами, ассоциированными
с топологиями Гротендика, возникают в контексте теорий
предпучкового типа. Фактически, если $\mathbb T$ --- теория
предпучкового типа, то теорема дуальности поднимается до биекции
между частными (их классами $\equiv_{s}$-эквивалентности) для
$\mathbb T$ и топологиями Гротендика $J$ на категории, обратной к
категории $\textrm{f.p.} {\mathbb T}\textrm{-mod}(\Set)$ конечно
представимых $\mathbb T$-моделей. Оказывается, что эта биекция
может быть сделана в теоретико-доказательной эквивалентности
способом, похожим на вышеописанный. Этот факт имеет важные
последствия; например, из него следует, что любое частное теории
предпучкового типа $\mathbb T$ имеет представление над своей
сигнатурой, в котором все аксиомы имеют вид ${\phi
\vdash_{\vec{x}} \mathbin{\mathop{\textrm{\huge
$\vee$}}\limits_{i\in I}}(\exists \vec{y_{i}})\theta_{i}}$, где
для любого ${i\in I}$ $\theta_{i}(\vec{y_{i}}, \vec{x})$ ---
$\mathbb T$-доказуемо функциональная формула из ${\{\vec{y_{i}}.
\psi\}}$ в ${\{\vec{x}. \phi\}}$, и $\phi(\vec{x})$,
$\psi(\vec{y_{i}})$ --- формулы, которые представляют $\mathbb
T$-модель.

Эти `теоретико-доказательные эквивалентности' также очень полезны
на практике, потому что они позволяют работать с топологиями
Гротендика в порядке получения аксиоматизаций геометрических
теорий или синтаксических результатов о них, и топологии
Гротендика имеют во многих контекстах замечательное
`вычислительное преимущество' над классическими системами
доказательств. Например, имеется полезная формула для топологии
Гротендика, образованной заданным семейством решет (это, например,
было использовано в \cite{OC} для установления теоремы дедукции
для геометрической логики) и для операции Гейтинга между
топологиями Гротендика (см. в \cite{OC} главу~3, содержащую все
эти результаты). Приложение этого дискурса к вычислению
(теоретико-решеточного) пересечения геометрических теорий было
дано в \cite{OC}, где естественная аксиоматизация пересечения
теории локальных колец и интегральных областей достигается
вычислением пересечения соответствующих топологий Гротендика на
категории, обратной к категории конечно представленных
коммутативных колец с единицей (заметим, что этот вид задач
аксиоматизации, вообще говоря, далек от тривиальности). Фактически
этот пример представляет собой приложение нашей философии `топосы
как мосты', описанной в разделе \ref{bridge}; действительно, в
этом случае инвариант есть пересечение двух подтопосов, которые с
точки зрения синтаксического сайта теории $\mathbb R$
коммутативных колец с единицей составляют пересечение двух частных
теории $\mathbb R$, тогда как с точки зрения тривиального сайта
классифицирующего топоса для $\mathbb R$ (категория $\cal R$,
обратная категории конечно представленных колец, оснащенных
тривиальной топологией Гротендика) этот инвариант соответствует
пересечению ассоциированных топологий Гротендика на категории
$\cal R$.

Продолжим описывать другие приложения теоремы дуальности. Из
теоремы следует, что набор (классов $\equiv_{s}$-эквивалентности)
геометрических теорий над данной сигнатурой имеет структуру
алгебры Гейтинга по отношению к естественному упорядочиванию
теорий, имеющему вид `${{\mathbb T}'\leq {\mathbb T}''}$, если и
только если все аксиомы из ${\mathbb T}'$ доказуемы в ${\mathbb
T}''$'. Также возможно вывести из теоремы явные описания решетки и
операций Гейтинга на геометрических теориях (см. \cite{OC}),
используя явные описания операций Гейтинга на топологиях
Гротендика. Это другая иллюстрация того факта, что работа с
топологиями Гротендика может быть легче, чем аргументация в
классической системе доказательств геометрической логики;
действительно, было бы весьма затруднительно достигнуть этих
характеризаций средствами традиционных логических методов.

В свете того факта, что понятие подтопоса есть теоретико-топосный
инвариант, теорема дуальности позволяет нам легко переносить
информацию между частными геометрических теорий, классифицируемых
одним и тем же топосом. Например, рассмотрим следующую задачу.
Пусть ${\mathbb T}_{1}$ и ${\mathbb T}_{2}$ --- две геометрические
теории, снабженные Морита-эквивалентностью между ними; верно ли,
что для любого частного ${\mathbb S}_{1}$ для ${\mathbb T}_{1}$
существует частное ${\mathbb S}_{2}$ для ${\mathbb T}_{2}$ такое,
что Морита-эквивалентность между ${\mathbb T}_{1}$ и ${\mathbb
T}_{2}$ ограничивается до Морита-эквивалентности между ${\mathbb
S}_{1}$ и ${\mathbb S}_{2}$? Теорема дуальности дает прямой
положительный ответ на этот вопрос. Действительно, подтопосы
классифицирующего топоса для ${\mathbb T}_{1}$ по теореме
дуальности соответствуют, с одной стороны, частным для ${\mathbb
T}_{1}$, а с другой стороны, поскольку классифицирующий топос для
${\mathbb T}_{1}$ эквивалентен классифицирующему топосу для
${\mathbb T}_{2}$, частным для ${\mathbb T}_{2}$. Отметим роль
классифицирующего топоса как моста (где инвариантом в данном
случае является понятие подтопоса, и два различных сайта
определения классифицирующего топоса являются геометрическими
синтаксическими сайтами этих двух теорий). Понимание, привнесенное
теоремой дуальности в этот контекст, не является только
теоретическим; как мы уже заметили, биекция между частными и
подтопосами часто может использоваться для получения явных
аксиоматизаций для геометрических теорий; в частности, явная
аксиоматизация для ${\mathbb S}_{2}$ начинается с аксиоматизации
для ${\mathbb S}_{1}$.

Вернемся теперь к фундаментальной роли теоремы дуальности в
вопросе возможности переноса концепций и результатов из
элементарной теории топосов в геометрическую логику.

Открытые (соответственно, замкнутые) подтопосы по дуальности
соответствуют частным, получаемым добавлением секвенций вида
${\top \vdash_{[]} \phi}$ (соответственно, ${\phi \vdash_{[]}
\bot}$), где $\phi$ --- геометрическое высказывание над сигнатурой
теории, и факторизация сюръекция-включение геометрического
морфизма имеет следующую естественную семантическую
интерпретацию.\\

\noindent
 {\bf Теорема~8.4.}
 {\it
 Пусть $\mathbb T$ --- геометрическая теория над сигнатурой
$\Sigma$, $f:{\cal F}\to {\cal E}$ --- геометрический морфизм в
классифицирующий топос $\cal E$ для $\mathbb T$, соответствующий
$\mathbb T$-модели $M$ в $\cal F$, как выше. Тогда топос ${\cal
E}'$ в факторизации сюръекция-включение ${\cal F}\epi {\cal E}'
\hookrightarrow {\cal E}$ для $f$ соответствует по теореме
дуальности частному $Th(M)$ для $\mathbb T$, состоящему из всех
геометрических секвенций $\sigma$ над $\Sigma$, которые имеются в
$M$.
 }\\

В частности, интересным теоретико-топосным инвариантом, который
посредством теоремы дуальности приводит к важным глубоким
результатам в геометрической теории, является \emph{булеанизация}
топоса, т.е. подтопос $\neg\neg\,$-пучков. Для данной
геометрической теории $\mathbb T$ назовем частное для $\mathbb T$,
соответствующее по теореме дуальности булеанизации
классифицирующего топоса для $\mathbb T$, \emph{булеанизацией}
$\mathbb T$.

Следующий результат (\cite{OC}, теорема 6.4.7) дает явную
аксиоматизацию для булеанизации геометрической теории.\\

\noindent
 {\bf Теорема~8.5.}
 {\it
 Пусть $\mathbb T$ --- геометрическая теория над сигнатурой
$\Sigma$. Тогда булеанизация $\mathbb T$ --- это теория,
полученная добавлением к аксиомам из $\mathbb T$ всех
геометрических секвенций вида ${\top \vdash_{\vec{x}}
\phi(\vec{x})}$, где $\phi(\vec{x})$ --- геометрическая формула
над $\Sigma$ такая, что для любой геометрической формулы
$\chi(\vec{x})$ над $\Sigma$ из доказуемости ${\phi \wedge \chi
\vdash_{\vec{x}} \bot}$ в $\mathbb T$ следует доказуемость ${\chi
\vdash_{\vec{x}} \bot}$ в $\mathbb T$.
 }\\

В случае важных математических теорий булеанизация часто приводит
к интересным частным. Например, булеанизация теории линейных
порядков (linear orders) --- это теория плотных линейных порядков
без конечных точек (endpoints), булеанизация теории булевых алгебр
--- это теория безатомных (atomless) булевых алгебр
(см. главу 9 из \cite{OC}), булеанизация (когерентной) теории
полей --- это (геометрическая) теория алгебраически замкнутых
полей конечной характеристики, в которых каждый элемент является
алгебраическим над первичным полем (см. \cite{OC4}). Более того,
при подходящих предположениях они аксиоматизируют слабо однородные
модели в смысле классической теории моделей (см. ниже разделы
\ref{presheaf} и \ref{fraisse}).

Другой теоретико-топосный инвариант, который, как уже доказано,
уместен в классической математике, --- это \emph{деморганизация}
топоса. Этот инвариант был введен в \cite{OC}, где было показано,
что каждый элементарный топос имеет наибольший плотный подтопос,
удовлетворяющий закону де Моргана; в этом же контексте название
\emph{деморганизация} было использовано для обозначения этого
подтопоса. В главе 6 из \cite{OC} была получена явная
аксиоматизация для частного геометрической теории,
соответствующего --- по теореме дуальности --- деморганизации
классифицирующего топоса этой теории, а в \cite{OC4} авторы
доказали, что деморганизация (когерентной) теории полей есть
геометрическая теория полей конечной характеристики, в которых
каждый элемент алгебраичен над первичным полем.

Данные два примера показывают, что эти мотивируемые логически
инварианты топосов имеют интересные проявления во множестве
математических контекстов; то, что они приводят к важным
математическим теориям, является ясным свидетельством их
центральной роли в математике.

\section{Примеры}

В этом разделе мы обсуждаем дальнейшие примеры из \cite{OC},
иллюстрирующие применение нашей философии `топосы как мосты'.

\subsection{Теории предпучкового типа}\label{presheaf}

Теория имеет \emph{предпучковый тип}, если это геометрическая
теория, классифицирующий топос которой эквивалентен предпучковому
топосу. Другими словами, геометрическая теория имеет предпучковый
тип, если и только если она Морита-эквивалентна теории плоских
функторов на малой категории $\cal C$. Класс теорий предпучкового
типа интересен по многим причинам. Одна важная концептуальная
причина состоит в том, что любая малая категория $\cal C$ может
рассматриваться --- с точностью до пополнения Коши
(Cauchy-completion) --- как категория $\textrm{f.p.} {\mathbb
T}\textrm{-mod}(\Set)$ (представлений классов изоморфизмов)
конечно представимых моделей теории предпучкового типа $\mathbb T$
(возьмем в качестве $\mathbb T$ теорию плоских функторов на $\cal
C$); как следствие, любой топос Гротендика имеет вид (с точностью
до эквивалентности) $\Sh(\textrm{f.p.} {\mathbb
T}\textrm{-mod}(\Set)^{\textrm{op}}, J)$ для некоторой теории
предпучкового типа $\mathbb T$. Другая причина в том, что этот
класс теорий содержит все декартовы, и в частности все конечно
алгебраические теории, а также много других интересных
математических теорий (например, теория линейных порядков, см.
\cite{MM} и \cite{El}, и геометрическая теория конечных множеств,
см. \cite{El}).

Классифицирующий топос теории предпучкового типа $\mathbb T$
всегда задается категорией функторов $[\textrm{f.p.} {\mathbb
T}\textrm{-mod}(\Set), \Set]$, где $\textrm{f.p.} {\mathbb
T}\textrm{-mod}(\Set)$ --- категория конечно представимых моделей
для $\mathbb T$ (см. \cite{OC}, глава 4). Следовательно, по
теореме дуальности частные для $\mathbb T$ биективно соответствуют
топологиям Гротендика на категории, обратной категории
$\textrm{f.p.} {\mathbb T}\textrm{-mod}(\Set)$; мы будем ссылаться
на топологию Гротендика, соответствующую частному ${\mathbb T}'$,
как на топологию Гротендика, \emph{ассоциированную} с ${\mathbb
T}'$. Это естественно приводит к следующему вопросу: могут ли
модели частного ${\mathbb T}'$ для $\mathbb T$ (в любом топосе
Гротендика) среди моделей для $\mathbb T$ быть характеризуемы
непосредственно в терминах ассоциированной топологии Гротендика?
Положительный ответ на этот вопрос дан в \cite{OC}, глава~4; а
именно, введено понятие \emph{$J$-однородной} модели теории
предпучкового типа $\mathbb T$ в топосе Гротендика (для топологии
Гротендика $J$ на $\textrm{f.p.} {\mathbb
T}\textrm{-mod}(\Set)^{\textrm{op}}$) и показано, что в каждом
топосе Гротендика модели частного ${\mathbb T}'$ для $\mathbb T$
--- это в точности $J$-однородные модели, где $J$ --- топология
Гротендика, ассоциированная с ${\mathbb T}'$. Если категория
$\textrm{f.p.} {\mathbb T}\textrm{-mod}(\Set)^{\textrm{op}}$
удовлетворяет правому условию Ора, и $J_{at}$ --- атомная
топология на ней, то понятие $J_{at}$-однородной модели
специализируется в $\Set$ понятием (слабо) однородной модели в
классической теории моделей (см. раздел~\ref{fraisse}).

Исследование класса теорий предпучкового типа дало автору работы
\cite{OC} прекрасную возможность протестировать методологии,
основанные на ее точке зрения `топосы как мосты', описанной выше.
Действительно, понятие теории предпучкового типа автоматически
поднимается до Морита-эквивалентности; а именно, если $\mathbb T$
--- теория предпучкового типа, то мы имеем Морита-эквивалентность
между $\mathbb T$ и теорией плоских функторов на категории
$\textrm{f.p.} {\mathbb T}\textrm{-mod}(\Set)^{\textrm{op}}$. В
теоретико-топосных терминах мы имеем эквивалентность топосов
Гротендика ${\Sh({\cal C}_{\mathbb T}, J_{\mathbb T}) \simeq
[\textrm{f.p.} {\mathbb T}\textrm{-mod}(\Set), \Set]}$. Два сайта
определения для этого классифицирующего топоса `достаточно
различны' для нетривиального переноса информации от одного к
другому, достижимого с использованием наших методов. Фактически в
\cite{OC} получено много результатов такого типа, кульминация
--- теоретико-топосная интерпретация
конструкции Фрессе в теории моделей (см. раздел \ref{fraisse}).
Сейчас мы коротко рассмотрим некоторые из них, чтобы показать, как
они возникают из применения наших принципов.

Мы уже обсуждали использование двойного представления
\[
\Sh({\cal C}_{\mathbb T}, J_{\mathbb T}) \simeq [\textrm{f.p.} {\mathbb T}\textrm{-mod}(\Set), \Set]
\]
классифицирующего топоса для $\mathbb T$ с целью аксиоматизации
частных для~$\mathbb T$.

Следующий результат дает связь между синтаксическими свойствами
частного ${\mathbb T}'$ для $\mathbb T$ и `топологическими'
свойствами ассоциированной топологии Гротендика $J$ на
$\textrm{f.p.} {\mathbb T}\textrm{-mod}(\Set)^{\textrm{op}}$.\\

\noindent
 {\bf Теорема~9.1.}
 {\it
 Пусть $\mathbb T$ --- теория предпучкового типа над сигнатурой
$\Sigma$, ${\mathbb T}'$ --- частное для $\mathbb T$ с
ассоциированной топологией Гротендика $J$ на $\textrm{f.p.}
{\mathbb T}\textrm{-mod}(\Set)^{\textrm{op}}$, и $\phi(\vec{x})$
--- геометрическая формула над $\Sigma$, которая представляет $\mathbb
T$-модель $M$. Тогда

\begin{enumerate}[(i)]

\item $\phi(\vec{x})$ $\mathbb T$-несводима; в частности,
$\phi(\vec{x})$ $\mathbb T$-доказуемо эквивалентна регулярной
формуле;

\item если сайт $(\textrm{f.p.} {\mathbb
T}\textrm{-mod}(\Set)^{\textrm{op}}, J)$ локально связен
(например, если $\textrm{f.p.} {\mathbb
T}\textrm{-mod}(\Set)^{\textrm{op}}$ удовлетворяет правому условию
Ора и каждое $J$-накрывающее решето непусто), то $\phi(\vec{x})$
${\mathbb T}'$-неразложима;

\item если $(\textrm{f.p.} {\mathbb
T}\textrm{-mod}(\Set)^{\textrm{op}}$ удовлетворяет правому условию
Ора и $J$ --- атомная топология на $(\textrm{f.p.} {\mathbb
T}\textrm{-mod}(\Set)^{\textrm{op}}$, то $\phi(\vec{x})$ ${\mathbb
T}'$-полна;

\item если каждое $J$-накрывающее решето на $M$ содержит
${J\mbox{-накрывающее}}$ решето, порожденное конечным семейством
морфизмов, то $\phi(\vec{x})$ ${\mathbb T}'$-компактна.

\end{enumerate}
 }

Доказательство этих результатов состоит в переносе (свойств)
некоего инварианта, а именно интерпретации формулы $\phi(\vec{x})$
в универсальной модели для ${\mathbb T}'$, сквозь два различных
представления классифицирующего топоса для ${\mathbb T}'$, таких
как категория $\Sh(\textrm{f.p.} {\mathbb
T}\textrm{-mod}(\Set)^{\textrm{op}}, J)$ и категория пучков
$\Sh({\cal C}_{{\mathbb T}'}, J_{{\mathbb T}'})$ на геометрическом
синтаксическом сайте для ${\mathbb T}'$. Конкретно, если
\[
l^{\textrm{f.p.} {\mathbb T}\textrm{-mod}(\Set)^{\textrm{op}}}_{J}:\textrm{f.p.} {\mathbb T}\textrm{-mod}(\Set)\to \Sh(\textrm{f.p.} {\mathbb T}\textrm{-mod}(\Set)^{\textrm{op}}, J)
\]
есть композиция функтора ассоциированного пучка
\[
[\textrm{f.p.} {\mathbb T}\textrm{-mod}(\Set), \Set] \rightarrow \Sh(\textrm{f.p.} {\mathbb T}\textrm{-mod}(\Set)^{\textrm{op}}, J)
\]
с вложением Йонеды
\[
\textrm{f.p.} {\mathbb T}\textrm{-mod}(\Set)^{\textrm{op}} \rightarrow [\textrm{f.p.} {\mathbb T}\textrm{-mod}(\Set), \Set]
\]
и
\[
y:{\cal C}_{{\mathbb T}'}\hookrightarrow \Sh({\cal C}_{{\mathbb T}'}, J_{{\mathbb T}'})
\]
--- вложение Йонеды, то имеется эквивалентность
\[
\Sh({\cal C}_{{\mathbb T}'}, J_{{\mathbb T}'}) \simeq
\Sh(\textrm{f.p.} {\mathbb T}\textrm{-mod}(\Set)^{\textrm{op}},
J)\, ,
\]
которая переводит ${y(\{\vec{x}. \phi\})}$ в $l^{\textrm{f.p.}
{\mathbb T}\textrm{-mod}(\Set)^{\textrm{op}}}_{J}(M)$, где $M$ ---
${\mathbb T\mbox{-модель}}$, представленная формулой
$\phi(\vec{x})$. Далее, перенос имеет место таким образом:
свойства сайта $(\textrm{f.p.} {\mathbb
T}\textrm{-mod}(\Set)^{\textrm{op}}, J)$ влекут за собой свойства
объекта $l^{\textrm{f.p.} {\mathbb
T}\textrm{-mod}(\Set)^{\textrm{op}}}_{J}(M)$ в классифицирующем
топосе $\Sh(\textrm{f.p.} {\mathbb
T}\textrm{-mod}(\Set)^{\textrm{op}}, J)$, которые переводятся с
помощью эквивалентности ${\Sh({\cal C}_{{\mathbb T}'}, J_{{\mathbb
T}'}) \simeq \Sh(\textrm{f.p.} {\mathbb
T}\textrm{-mod}(\Set)^{\textrm{op}}, J)}$ в свойства объекта
${y(\{\vec{x}. \phi\})}$ в классифицирующем топосе $\Sh({\cal
C}_{{\mathbb T}'}, J_{{\mathbb T}'})$, которые в свою очередь
переводятся в свойства синтаксического сайта $({\cal C}_{{\mathbb
T}'}, J_{{\mathbb T}'})$ и, следовательно, теории ${\mathbb T}'$.

Следующий результат, выражающий тот факт, что синтаксические и
семантические понятия конечной представимости модели (формальное
определение этих концепций см. в \cite{OC}) совпадают для теорий
предпучкового типа, получен тем же самым методом с использованием
инвариантного свойства объекта из топоса Гротендика быть
неприводимым.\\

\noindent
 {\bf Теорема~9.2.}
 {\it
 Пусть $\mathbb T$ --- теория предпучкового типа над сигнатурой
$\Sigma$. Тогда

\begin{enumerate}

\item любая конечно представимая $\mathbb T$-модель в $\Set$
представляется ${\mathbb T\mbox{-несводимой}}$ геометрической
формулой $\phi(\vec{x})$ над $\Sigma$;

\item обратно, любая $\mathbb T$-несводимая геометрическая формула
$\phi(\vec{x})$ над $\Sigma$ представляет $\mathbb T$-модель.
\end{enumerate}
В частности, категория $\textrm{f.p.} {\mathbb
T}\textrm{-mod}(\Set)^{\textrm{op}}$ эквивалентна полной
подкатегории из ${\cal C}_{\mathbb T}^{\textrm{geom}}$ на $\mathbb
T$-несводимых формулах.
 }\\

Другой интересный результат о теориях предпучкового типа
--- следующая теорема определимости (\cite{OC},
следствие~7.2.2).\\

\noindent
 {\bf Теорема~9.3.}
 {\it
 Пусть $\mathbb T$ --- теория предпучкового типа над сигнатурой
$\Sigma$, пусть $A_{1}, \ldots, A_{n}$ --- строка типов для
$\Sigma$, и пусть для любой конечно представимой $\Set$-модели $M$
для $\mathbb T$ задано подмножество $R_{M}$ из ${MA_{1}\times
\cdots \times MA_{n}}$ таким способом, что каждый гомоморфизм
${\mathbb T\mbox{-моделей}}$ ${h:M\to N}$ отображает $R_{M}$ в
$R_{N}$. Тогда существует геометрическая формула-в-контексте
${\phi(x^{A_{1}}, \ldots, x^{A_{n}})}$ такая, что
${R_{M}=[[\phi]]_{M}}$ для любой $M$.
 }\\

Этот результат опять доказывается методом переноса инварианта
сквозь различные сайты определения данного классифицирующего
топоса. В этом случае инвариантом является понятие подобъекта
рассматриваемого объекта (объектов) универсальной модели~$\mathbb
T$, и два представления классифицирующего топоса для $\mathbb T$
--- это $[\textrm{f.p.} {\mathbb
T}\textrm{-mod}(\Set)^{\textrm{op}}, \Set]$ и $\Sh({\cal
C}_{\mathbb T}, J_{\mathbb T})$. Аргументация следующая. В
терминах представления $[\textrm{f.p.} {\mathbb
T}\textrm{-mod}(\Set)^{\textrm{op}}, \Set]$ универсальная модель
для $\mathbb T$ может быть описана как функтор $M_{\mathbb T}$ в
$[\textrm{f.p.} {\mathbb T}\textrm{-mod}(\Set), \Set]$, который
сопоставляет типу $A$ функтор $M_{\mathbb T}A$ вида $(M_{\mathbb
T}A)(M)=MA$, определяемый для функции и символов отношения
очевидным способом, тогда как в терминах другого представления
$\Sh({\cal C}_{\mathbb T}, J_{\mathbb T})$ универсальная модель
для $\mathbb T$ допускает такое синтаксическое описание, что
интерпретация формулы $\phi(\vec{x})$ в ней имеет вид
${y(\{\vec{x}. \phi\})}$ (см. детали в \cite{OC}, глава~7). Тогда
теорема следует как раз из замечания, что сопоставление ${M\to
R_{M}}$ в утверждении теоремы поднимается до подобъекта ${R\mono
M_{\mathbb T}}$ в топосе $[\textrm{f.p.} {\mathbb
T}\textrm{-mod}(\Set), \Set]$.

Отметим, что этот результат тесно связан с особой природой теорий
предпучкового типа. Можно ожидать, что та же самая теорема будет
справедлива для б\'{о}льших классов теорий; например, теорема,
вообще говоря, не имеет места для когерентных теорий.
Действительно, теория разрешимых (decidable) колец (т.е. теория
коммутативных колец с единицей, снабженных предикатом, который
доказуемо дополняется до отношения равенства), очевидно,
удовлетворяет предположениям теоремы как для свойства элемента
быть нильпотентом, так и для дополнения к этому свойству, и если
последнее свойство было определимо с помощью геометрической
формулы, то каждое свойство должно было бы быть определимо при
помощи когерентной формулы, а это, как может быть легко доказано
(в смысле компактности аргументации), неверно.

Мы отметим, что эти результаты, которые, как мы видели,
естественно следуют из реализации нашей идеи `топосы как мосты',
являются, несомненно, тривиальными, тогда как, будучи
специализируемы полезными результатами во множестве математических
контекстов, они могли бы быть очень трудны для доказательства с
использованием традиционной техники. С другой стороны, бесконечное
число результатов может быть порождено почти `автоматически' с
использованием тех же самых принципов для различных топосов и
различных инвариантов; количество нетривиальных математических
результатов, которые могут быть `автоматически генерированы' с
применением этой методологии, просто безгранично. Эти методы
представляют собой новый способ порождения математики, который как
формально, так и по существу отличается от того, что дает
классическая техника.

\subsection{Конструкция Фрессе в теоретико-топосной перспективе}\label{fraisse}

В этом разделе мы обсуждаем особо значимый пример (из \cite{OC},
глава~9) результата, который, хотя и был получен применением наших
чисто теоретико-топосных методов, прямо соотносится с хорошо
известным разделом математики, а именно с конструкцией Фресс\'{e}
(Fra\"{\i}ss\'e) в теории моделей.

В разделе \ref{presheaf} мы обсуждали приложение философии `топосы
как мосты' в контексте инвариантов \emph{объектов} топосов; вместо
этого главный результат настоящего раздела возникает из
рассмотрения `глобальных' инвариантов на данном классифицирующем
топосе.

Прежде чем сформулировать теорему, дадим несколько определений.\\

\noindent
 {\bf Определение~9.4.}
 Категория $\cal C$ называется удовлетворяющей \emph{свойству
амальгамы}, если для любых объектов $a,b,c\in {\cal C}$ и
морфизмов $f:a\rightarrow b$, $g:a\rightarrow c$ из $\cal C$
существуют объект $d\in \cal C$ и морфизмы $f':b\rightarrow d$,
$g':c\rightarrow d$ в $\cal C$ такие, что $f'\circ f=g'\circ g$:
\[
\xymatrix {
a \ar[d]_{g} \ar[r]^{f} & b  \ar@{-->}[d]^{f'} \\
c \ar@{-->}[r]_{g'} & d }
\]

Отметим, что $\cal C$ удовлетворяет свойству амальгамы, если и
только если ${\cal C}^{\textrm{op}}$ удовлетворяет правому условию
Ора. Поэтому если $\cal C$ удовлетворяет свойству амальгамы, то мы
можем снабдить ${\cal C}^{\textrm{op}}$ атомной топологией.\\

\noindent
 {\bf Определение~9.5.}
 Категория $\cal C$ называется удовлетворяющей \emph{свойству
совместного погружения}, если для любой пары объектов ${a,b\in
{\cal C}}$ существуют объект ${c\in \cal C}$ и морфизмы
${f:a\rightarrow c}$, ${g:b\rightarrow c}$ из $\cal C$:
\[
\xymatrix {
 & a  \ar@{-->}[d]^{f} \\
b \ar@{-->}[r]_{g} & c }
\]

\noindent
 {\bf Определение~9.6.}
 Пусть $\mathbb T$ --- теория предпучкового типа. Модель $M$ для
$\mathbb T$ называется \emph{однородной}, если для любых моделей
${a,b \in \textsf{f.p.} {\mathbb T}\textsf{-mod}(\Set)}$ и стрелок
${j:a\rightarrow b}$ и ${\chi:a\rightarrow M}$ из ${\mathbb
T}\textsf{-mod}(\Set)$ существует стрелка
${\tilde{\chi}:b\rightarrow M}$ из ${\mathbb
T}\textsf{-mod}(\Set)$ такая, что $\tilde{\chi}\circ j=\chi$:
\[
\xymatrix {
a \ar[d]_{j} \ar[r]^{\chi} & M \\
b \ar@{-->}[ur]_{\tilde{\chi}} &  }
\]

Наш главный результат следующий.\\

\noindent
 {\bf Теорема~9.7.}
 {\it
 Пусть $\mathbb{T}$ --- теория предпучкового типа такая,
что категория $\textrm{f.p.} {\mathbb T}\textrm{-mod}(\Set)$
удовлетворяет свойствам амальгамы и совместного погружения. Тогда
любые две вычислимые однородные ${\mathbb T\mbox{-модели}}$ в
$\Set$ изоморфны.
 }\\

Мы можем описать структуру доказательства этой теоремы следующим
образом.

Если категория $\textrm{f.p.} {\mathbb T}\textrm{-mod}(\Set)$
пуста, то $\mathbb T$ не имеет моделей в $\Set$, и, значит,
утверждение тривиально выполняется. Поэтому мы будем предполагать,
что $\textrm{f.p.} {\mathbb T}\textrm{-mod}(\Set)$ непуста. Так
как $\textrm{f.p.} {\mathbb T}\textrm{-mod}(\Set)$ удовлетворяет
свойству амальгамы, мы можем снабдить обратную категорию
$\textrm{f.p.} {\mathbb T}\textrm{-mod}(\Set)^{\textrm{op}}$
атомной топологией $J_{at}$. Наша аргументация основана на топосе
$\Sh(\textrm{f.p.} {\mathbb T}\textrm{-mod}(\Set)^{\textrm{op}},
J_{at})$. По теореме дуальности подтопос
\[
\Sh(\textrm{f.p.} {\mathbb T}\textrm{-mod}(\Set)^{\textrm{op}}, J_{at}) \hookrightarrow [\textrm{f.p.} {\mathbb T}\textrm{-mod}(\Set)^{\textrm{op}}, \Set]
\]
соответствует единственному частному (с точностью до
синтаксической эквивалентности) ${\mathbb T}'$ для $\mathbb T$
такому, что его классифицирующим топосом является
$\Sh(\textrm{f.p.} {\mathbb T}\textrm{-mod}(\Set)^{\textrm{op}},
J_{at})$ (Заметим, что концептуально использование теоремы
дуальности существенно состоит в переносе инвариантного понятия
подтопоса из представления $[\textrm{f.p.} {\mathbb
T}\textrm{-mod}(\Set)^{\textrm{op}}, \Set]$ классифицирующего
топоса для $\mathbb T$ в синтаксическое его представление
$\Sh({\cal C}_{\mathbb T}, J_{\mathbb T})$). Из техники
`представления Йонеды для плоских функторов' (см. \cite{OC},
глава~4) следует, что в каждом топосе Гротендика, в частности в
$\Set$, модели для ${\mathbb T}'$ --- это в точности однородные
$\mathbb T$-модели.

Далее, `ядро' аргументации состоит в переносе инвариантного
свойства топоса быть \emph{атомным и двузначным} топосом из
`семантического' представления $\Sh(\textrm{f.p.} {\mathbb
T}\textrm{-mod}(\Set)^{\textrm{op}}, J_{at})$ в синтаксическое
представление $\Sh({\cal C}_{{\mathbb T}'}, J_{{\mathbb T}'})$
классифицирующего топоса для ${\mathbb T}'$. Можно доказать, что
если $\textrm{f.p.} {\mathbb T}\textrm{-mod}(\Set)$ удовлетворяет
свойству совместного погружения, равно как и свойству амальгамы,
тогда топос $\Sh(\textrm{f.p.} {\mathbb
T}\textrm{-mod}(\Set)^{\textrm{op}}, J_{at})$ двузначен, так же
как и является атомным (фактически второе утверждение тоже имеет
место при наших предположениях). Но если прочесть этот инвариант с
точки зрения другого представления $\Sh({\cal C}_{{\mathbb T}'},
J_{{\mathbb T}'})$, это в точности скажет нам, что теория
${\mathbb T}'$ `атомна' и `полна', что обязывает ${\mathbb T}'$
быть вычислимо категорной (см. \cite{OC}, глава~8). Если принять
во внимание вышеупомянутый факт, что модели для ${\mathbb T}'$ ---
в точности однородные $\mathbb T$-модели, это завершает
доказательство теоремы. Из переноса инвариантной конструкции
булеанизации сквозь два различных представления классифицирующего
топоса для $\mathbb T$ сразу же видно, что частное ${\mathbb T}'$
--- это в точности булеанизация для $\mathbb T$.

Другое интересное применение нашей методологии переноса
инвариантов сквозь различные сайты определения одного и того же
топоса (данной в \cite{OC}, глава~9) касается существования
однородных моделей теории предпучкового типа $\mathbb T$ в $\Set$.
Как мы уже отмечали, имеется частное ${\mathbb T}'$ для $\mathbb
T$, которое аксиоматизирует однородные модели в каждом топосе
Гротендика и классифицируется топосом $\Sh(\textrm{f.p.} {\mathbb
T}\textrm{-mod}(\Set)^{\textrm{op}}, J_{at})$. Следовательно,
существование однородной модели для $\mathbb T$ в $\Set$
эквивалентно существованию точки топоса $\Sh(\textrm{f.p.}
{\mathbb T}\textrm{-mod}(\Set)^{\textrm{op}}, J_{at})$. Заметим,
что это последнее свойство является теоретико-топосным
инвариантом. Далее, можно легко доказать, что если топос $\cal E$
имеет достаточно точек, то $\cal E$ имеет точку тогда и только
тогда, когда он нетривиален; и можно непосредственно доказать, что
топос $\Sh(\textrm{f.p.} {\mathbb
T}\textrm{-mod}(\Set)^{\textrm{op}}, J_{at})$ нетривиален при
условии, что категория $\textrm{f.p.} {\mathbb
T}\textrm{-mod}(\Set)$ непуста. Следовательно, при предположении,
что $\textrm{f.p.} {\mathbb T}\textrm{-mod}(\Set)$ непуста, все
сводится к доказательству, что $\Sh(\textrm{f.p.} {\mathbb
T}\textrm{-mod}(\Set)^{\textrm{op}}, J_{at})$ имеет достаточно
точек. Интересно, что это может быть сделано различными способами;
например, по теореме Делиня достаточно найти когерентную теорию,
которая классифицируется топосом $\Sh(\textrm{f.p.} {\mathbb
T}\textrm{-mod}(\Set)^{\textrm{op}}, J_{at})$. Фактически в
\cite{OC} мы пользовались тем, что чисто комбинаторное условие на
категорию $\textrm{f.p.} {\mathbb T}\textrm{-mod}(\Set)$
гарантирует, что теория $J_{at}$-непрерывных плоских функторов на
$\textrm{f.p.} {\mathbb T}\textrm{-mod}(\Set)^{\textrm{op}}$
когерентна.

\subsection{Другие примеры}\label{examp}

В этом разделе мы обсуждаем другие приложения философии `топосы
как мосты' к вопросам логики.

Напомним, что для любого фрагмента геометрической логики имеется
соответствующее понятие доказуемости для теорий в этом фрагменте;
например, имеется понятие доказуемости регулярных секвенций в
регулярной логике и понятие доказуемости когерентных секвенций в
когерентной логике, так же как и классическое понятие доказуемости
геометрических секвенций в геометрической логике.

Естественный вопрос: совместимы ли друг с другом эти понятия
доказуемости, т.е. приводит ли понятие доказуемости в данном
фрагменте логики к понятиям доказуемости в меньшем фрагменте
(заметим, что любая теория в данном фрагменте может представляться
как теория в большем фрагменте). Мы показываем, что имеется
естественный теоретико-топосный путь работы с этими проблемами,
который использует тот факт, что для теорий в подходящем фрагменте
геометрической логики имеется много синтаксических представлений
их классифицирующих топосов, каждое из которых соответствует
частному фрагменту логики, в котором лежит теория (см. выше
раздел~\ref{onetopos}). Например, если теория $\mathbb T$ над
сигнатурой $\Sigma$ когерентна, то имеются два синтаксических
представления ее классифицирующего топоса: как категории пучков
$\Sh({\cal C}_{\mathbb T}^{\textrm{coh}}, J_{{\cal C}_{\mathbb
T}^{\textrm{coh}}})$ на когерентном синтаксическом сайте $({\cal
C}_{\mathbb T}^{\textrm{coh}}, J_{{\cal C}_{\mathbb
T}^{\textrm{coh}}})$ для $\mathbb T$ и как категории пучков
$\Sh({\cal C}_{\mathbb T}, J_{{\cal C}_{\mathbb T}})$ на
геометрическом синтаксическом сайте $({\cal C}_{\mathbb T},
J_{{\cal C}_{\mathbb T}})$ для $\mathbb T$. Из того факта, что
вложения Йонеды ${y_{\textrm{coh}}:{\cal C}_{\mathbb
T}^{\textrm{coh}} \to \Sh({\cal C}_{\mathbb T}^{\textrm{coh}},
J_{{\cal C}_{\mathbb T}^{\textrm{coh}}})}$ и ${y:{\cal C}_{\mathbb
T}\to \Sh({\cal C}_{\mathbb T}, J_{{\cal C}_{\mathbb T}})}$
консервативны, следует, что когерентная секвенция над $\Sigma$
доказуема в универсальной модели для $\mathbb T$, лежащей в топосе
$\Sh({\cal C}_{\mathbb T}^{\textrm{coh}}, J_{{\cal C}_{\mathbb
T}^{\textrm{coh}}})$, если и только если она доказуема в $\mathbb
T$ с использованием когерентной логики, и что геометрическая
секвенция над $\Sigma$ доказуема в универсальной модели для
$\mathbb T$, лежащей в топосе $\Sh({\cal C}_{\mathbb
T}^{\textrm{coh}}, J_{{\cal C}_{\mathbb T}^{\textrm{coh}}})$, если
и только если она доказуема в $\mathbb T$ с использованием
геометрической логики. Но свойство действительности данного
секвента в универсальной модели геометрической теории очевидно
является теоретико-топосным инвариантом, откуда следует, что
когерентная секвенция доказуема в $\mathbb T$ с использованием
когерентной логики, если и только если она доказуема в $\mathbb T$
с использованием геометрической логики, что и требуется.

Другое интересное применение касается `классической полноты'
теорий по отношению к геометрической логике. Для когерентной
теории $\mathbb T$ над сигнатурой $\Sigma$ хорошо известно, что
если когерентная секвенция над $\Sigma$ удовлетворяется в каждой
основанной на $\Set$ модели для $\mathbb T$, то она доказуема в
$\mathbb T$ с использованием когерентной логики. Естественно,
таким образом, задаться вопросом: справедливо ли то же самое для
геометрических секвенций при том же самом предположении о
когерентности $\mathbb T$, т.е. доказуема ли в $\mathbb T$ с
использованием геометрической логики произвольная геометрическая
секвенция над $\Sigma$, которая действительна в любой основанной
на $\Set$ модели для $\mathbb T$. Ответ на этот вопрос
утвердителен и естественно следует из переноса инварианта, а
именно свойства топоса иметь достаточно точек, сквозь два
различных представления $\Sh({\cal C}_{\mathbb T}^{\textrm{coh}},
J_{{\cal C}_{\mathbb T}^{\textrm{coh}}}) \simeq \Sh({\cal
C}_{\mathbb T}, J_{{\cal C}_{\mathbb T}})$ классифицирующего
топоса для $\mathbb T$. Действительно, в терминах сайта $({\cal
C}_{\mathbb T}^{\textrm{coh}}, J_{{\cal C}_{\mathbb
T}^{\textrm{coh}}})$ инвариант перефразируется как условие, что
каждая когерентная секвенция над $\Sigma$, которая удовлетворяется
во всех основанных на $\Set$ моделях для $\mathbb T$, должна быть
доказуема в $\mathbb T$ с использованием когерентной логики, тогда
как в терминах сайта $({\cal C}_{\mathbb T}, J_{{\cal C}_{\mathbb
T}})$ он перефразируется как условие, что каждая геометрическая
секвенция над $\Sigma$, которая удовлетворяется во всех основанных
на $\Set$ моделях для $\mathbb T$, должна быть доказуема в
$\mathbb T$ с использованием геометрической логики.

Другой результат, вполне естественно вытекающий из двойного
представления $\Sh({\cal C}_{\mathbb T}^{\textrm{coh}}, J_{{\cal
C}_{\mathbb T}^{\textrm{coh}}}) \simeq \Sh({\cal C}_{\mathbb T},
J_{{\cal C}_{\mathbb T}})$ классифицирующего топоса когерентной
теории, представляет собой следующее (\cite{OC},
теорема~10.2.5).\\

 \noindent
 {\bf Теорема~9.8.}
 {\it
Пусть $\mathbb T$ --- геометрическая теория над сигнатурой
$\Sigma$. Тогда $\mathbb T$ когерентна, если и только если для
любой когерентной формулы ${\{\vec{x}. \phi\}}$ над $\Sigma$ и для
любого семейства ${\{\psi_{i}(\vec{x}) \textrm{ | } i\in I\}}$
геометрических формул в том же самом контексте из
$$
\phi \vdash_{\vec{x}} \mathbin{\mathop{\textrm{\huge
$\vee$}}\limits_{i\in I}} \psi_{i} \mbox{ доказуема в } \mathbb T
$$
следует
$$
\phi \vdash_{\vec{x}} \mathbin{\mathop{\textrm{\huge
$\vee$}}\limits_{i\in I'}} \psi_{i} \mbox{ доказуема в } \mathbb T
$$
для некоторого конечного подмножества $I'$ из $I$.
 }\\

Нетривиальная часть теоремы немедленно следует из замечания, что
если $\mathbb T$ --- когерентная теория над сигнатурой $\Sigma$,
то имеется эквивалентность классифицирующих топосов ${\Sh({\cal
C}_{\mathbb T}^{\textrm{coh}}, J_{{\cal C}_{\mathbb
T}^{\textrm{coh}}}) \simeq \Sh({\cal C}_{\mathbb T}, J_{{\cal
C}_{\mathbb T}})}$ такая, что для любой когерентной формулы
${\{\vec{x}. \phi\}}$ над $\Sigma$ объекты
${y_{\textrm{coh}}(\{\vec{x}.\phi\})}$ и ${y(\{\vec{x}.\phi\})}$
соответствуют друг другу при эквивалентности (здесь
${y_{\textrm{coh}}}$ и $y$ --- два вложения Йонеды
${y_{\textrm{coh}}:{\cal C}_{\mathbb T}^{\textrm{coh}} \to
\Sh({\cal C}_{\mathbb T}^{\textrm{coh}}, J_{{\cal C}_{\mathbb
T}^{\textrm{coh}}})}$ и ${y:{\cal C}_{\mathbb T}\to \Sh({\cal
C}_{\mathbb T}, J_{{\cal C}_{\mathbb T}})}$). Если мы рассмотрим
инвариантное свойство объекта из топоса быть компактным, мы
увидим, что из факта, что сайт $({\cal C}_{\mathbb
T}^{\textrm{coh}}, J_{{\cal C}_{\mathbb T}^{\textrm{coh}}})$
когерентен, следует, что объект
${y_{\textrm{coh}}(\{\vec{x}.\phi\})}$ компактен, но с точки
зрения геометрического синтаксического сайта $({\cal C}_{\mathbb
T}, J_{{\cal C}_{\mathbb T}})$ свойство объекта
${y(\{\vec{x}.\phi\})}$ быть компактным перефразируется так: для
любого семейства ${\{\psi_{i}(\vec{x}) \textrm{ | } i\in I\}}$
геометрических формул в том же самом контексте из
$$
\phi \vdash_{\vec{x}} \mathbin{\mathop{\textrm{\huge
$\vee$}}\limits_{i\in I}} \psi_{i} \mbox{ доказуема в } \mathbb T
$$
следует
$$
\phi \vdash_{\vec{x}} \mathbin{\mathop{\textrm{\huge
$\vee$}}\limits_{i\in I'}} \psi_{i} \mbox{ доказуема в } \mathbb T
$$
для некоторого конечного подмножества $I'$ из $I$.

Эта теорема показывает: то, что характеризует когерентные теории
среди геометрических, есть в точности общая форма компактности.

Конечно, все вышеупомянутые результаты имеют свою версию для
декартовых или регулярных теорий вместо когерентных (см.
\cite{OC}, глава~10).

Наконец, мы заметим, что хорошо известная эквивалентность между
теоремой Делиня в алгебраической геометрии, утверждающей, что
каждый когерентный топос имеет достаточно точек, и классической
теоремой о полноте для когерентных теорий в логике может быть
естественно выражена в наших рамках. Действительно, начиная с
когерентного топоса, представленного как категория пучков на сайте
$({\cal C}, J)$ таком, что $\cal C$ имеет конечные пределы, и $J$
--- топология Гротендика конечного типа на ней, можно рассмотреть
для него другой сайт определения, а именно когерентный
синтаксический сайт $({\cal C}_{\mathbb T}^{\textrm{coh}},
J_{{\cal C}_{\mathbb T}^{\textrm{coh}}})$ когерентной теории
$\mathbb T$ (так как согласно нашему предположению теория плоских
$J$-непрерывных функторов на $\cal C$ допускает когерентную
аксиоматизацию). Таким образом, для данного когерентного топоса
имеются два различных представления: одно --- `геометрической'
природы, и другое
--- логической природы; и
вышеупомянутая эквивалентность следует в точности из того факта,
что выражение инвариантного свойства иметь достаточно точек в
терминах когерентного синтаксического сайта $({\cal C}_{\mathbb
T}^{\textrm{coh}}, J_{{\cal C}_{\mathbb T}^{\textrm{coh}}})$ для
$\mathbb T$ составляет в точности классическую теорему полноты для
$\mathbb T$.

\section{Топосы для работающего математика}\label{work}

В этом разделе мы суммируем некоторые аспекты теории топосов, и в
частности --- той методологии, которую мы описали выше в данной
статье, чтобы сделать предмет привлекательным для `работающего
математика'.

Мы видели, что каждой математической теории первого порядка
(допускающей геометрическую аксиоматизацию) можно естественным
образом сопоставить классифицирующий топос, который выделяет
существенные детали теории, а именно те детали, которые
инвариантны по отношению к Морита-эквивалентности. Данный топос
может иметь много различных сайтов определения, которые
высвечивают различные аспекты терий, им классифицируемых.
Теоретико-топосные инварианты (т.е. свойства топосов или
конструкции, включающие их, которые инвариантны относительно
категорной эквивалентности) могут тогда быть использованы для
переноса свойств между теориями, классифицируемыми этим топосом.
Фактически абстрактное соотношение между сайтом $({\cal C},J)$ и
топосом $\Sh({\cal C}, J)$, который он `порождает', часто очень
естественно, позволяя нам легко переносить инварианты между
различными сайтами. Более того, как мы видели в разделах
\ref{logic} и \ref{duality}, уровень общности, представленный
теоретико-топосными инвариантами, идеален для отлавливания
некоторых важных деталей математических теорий.

Разновидность `кодирования' ${({\cal C}, J)\to \Sh({\cal C}, J)}$,
даваемая операцией взятия категории пучков на данном сайте, в
высшей степени нетривиальна; два ингредиента $({\cal C}, J)$
изощренным способом комбинируются вместе для создания топоса
$\Sh({\cal C}, J)$. Роль теоретико-топосных инвариантов состоит в
точности в `декодировании' информации, скрытой в топосе, в смысле
идентификации свойств сайта $({\cal C}, J)$, которые
предполагаются свойствами соответствующего топоса $\Sh({\cal C},
J)$, или описания результата построений, определенных на топосе,
непосредственно в терминах сайта способом, который понятен
специалистам вне теории топосов. Действительно, тогда как понятие
сайта вполне близко общей математической практике, топос является
в действительности объектом другой природы, особой сущностью,
которая может быть изучена посредством специфических методов (т.е.
методов теории топосов), которые \emph{не} являются методами
`классической математики' (см. раздел~\ref{genetics}, где имеется
концептуальная аналогия, которая может прояснить этот вопрос).

Связь между классической математикой и теорией топосов дается тем
фактом, что многие конструкции или свойства, естественно
возникающие в математике, могут быть выражены в терминах
теоретико-топосных инвариантов; например, группы когомологий
топосов являются теоретико-топосными инвариантами, которые в
случае конкретных топосов специализируются когомологиями,
возникающими в математической практике (такими как, например,
обычные когомологии топологических пространств, когомологии Галуа,
когомологии Эйленберга--Маклейна, когомологии Вейля), допуская в
то же время достаточно свободы для конструирования новых
`конкретных' когомологий, которые могут быть полезны для других
целей.

Далее, преимущество работы с топосами более чем с другими
разновидностями сущностей заключается в том факте, что --- как мы
видели в статье --- топосы могут естественно работать как
унифицирующие пространства, допуская перенос свойств или
конструкций между различными математическими теориями и открывая
таким образом путь к использованию методов одного раздела
математики для решения задач в другом. Например, свойство группы
когомологий данного топоса с успехом может быть переведено --- с
использованием представления топоса как категории пучков на
синтаксическом сайте геометрической теории $\mathbb T$ --- в
некоторое логическое свойство теории $\mathbb T$ (похожим образом
с использованием сайта алгебраической природы можно прийти к
свойству, трактуемому методами алгебры, и т.д.).

`Работающий математик' может очень успешно пытаться формулировать
интересующие его или ее свойства в терминах теоретико-топосных
инвариантов и порождать эквивалентные версии их с использованием
альтернативных сайтов. Также он может пробовать выразить
математическую двойственность или эквивалентность, которые он
хочет исследовать, в терминах Морита-эквивалентности, вычислить
классифицирующий топос двух теорий и применить теоретико-топосные
методы для извлечения новой информации о ней (см. разделы
\ref{onetopos}, \ref{morita} и \ref{bridge}).

Следовательно, может быть приложено много исследовательских усилий
для введения новых инвариантов топосов, которые интересно связаны
с вопросами, изучаемыми математиками, и идентификации
Морита-эквивалентностей (также в виде установления общих теорем
представления для топосов); мы надеемся убедить читателя, что
только лишь с использованием наиболее известных инвариантов
топосов можно легко получать существенные результаты в различных
математических контекстах. Можно искать для таких инвариантов
характеризации сайтов с целью сделать перенос свойств осуществимым
при наличии Морита-эквивалентности.

Мы отмечаем, что одна математическая теория может `порождать'
значительное число нетривиальных Морита-эквивалентностей (см.
раздел~\ref{onetopos}); таким образом, теоретико-топосные методы
могут быть использованы для извлечения информации с помощью только
одной теории без необходимости специфического поиска
Морита-эквивалентностей, связывающих эту теорию с `выглядящими
иначе' теориями.

Наконец, теория топосов действительно может предложить очень много
специалистам в любой области математики; в смысле методологий,
описанных в статье выше, можно порождать огромное число новых
результатов в любой области математики без каких-либо творческих
усилий. Этот унифицирующий механизм реально имеет потенциал для
автоматического генерирования результатов; конечно, многие
результаты, полученные таким способом, окажутся для работающего
математика `сверхъестественными' (хотя они могут быть еще и весьма
глубокими), и фактически с целью получения интересных результатов
нужно тщательно отбирать инварианты и Морита-эквивалентности так,
чтобы они естественно соотносились с интересующими вопросами (как,
например, в случае конструкции Фрессе). Еще надо иметь в виду, что
тот вид понимания, который эти методы могут принести в математику,
\emph{существенно} отличен от чего-либо, что можно было бы
естественно ожидать. В точности там, где традиционные методы
терпят неудачу, теоретико-топосные методы могут, возможно,
привести к решению (подумайте, например, о доказательстве
Гротендиком и Делинем гипотезы Вейля).

\subsection{Сравнение с генетикой}\label{genetics}

Та разновидность понимания, которую теоретико-топосные методы
могут привнести в математику, сравнима с тем, что генетика
привнесла в медицину (биологам приносятся извинения за возможные
неточности по поводу их предмета, могущие проявиться в следующем
абзаце).

Как (человеческая) ДНК содержит существенные детали индивидуума,
так и классифицирующий топос содержит существенные детали
математической теории. Как ДНК может быть выделена множеством
разных способов (например, из различных частей тела индивидуума),
так и классифицирующий топос может быть представлен и вычислен
альтернативными способами (выделяющими различные прочтения ядра
теории). Как ДНК инвариантна по отношению к частному
физиологическому проявлению индивидуума в данное время (например,
по отношению к возрасту), так и классифицирующий топос инвариантен
по отношению к отдельным представлениям теории (например, по
отношению к отдельной ее аксиоматизации над ее сигнатурой). Как
ДНК индивидуума может быть изучена с использованием подходящей
техники, не относящейся к традиционной медицине, так и
классифицирующий топос может быть изучен с использованием
специальных методов (т.е. методов теории топосов), которые, хоть и
являются полностью математическими, отличны от методов
классической математики. Как в генетике изучается влияние
изменений в ДНК на характеристики индивидуума, так и в теории
топосов изучается эффект, который теоретико-топосные операции на
топосах оказывают на теории, ими классифицируемые. Как роль ДНК
--- это роль унифицирующей концепции, позволяющей сравнивать
индивидуумов друг с другом, выделяя различия и вскрывая сходство,
так и понятие классифицирующего топоса является унифицирующим,
позволяя сравнивать различные математические теории друг с другом
и переносить знания между ними.

Мы надеемся передать этой метафорой интуицию, стоящую за
представленными результатами, а также потенциальную роль
теоретико-топосной техники в математике в будущем.

\vspace{0.5cm}
 {\flushleft
 {\bf Благодарности:}
 }
 Я признательна
Марко Бенини за его полезные комментарии к черновику этой статьи.

\end{document}